\documentclass[a4paper, leqno, 11pt]{amsart}
\usepackage[utf8]{inputenc}
\usepackage{amsmath, amsfonts, amssymb}
\usepackage{geometry}
\usepackage{graphicx}
\usepackage{enumerate}

\usepackage{todonotes}

\newcommand{\halb}{\frac{1}{2}}
\theoremstyle{definition}
\newtheorem{algo}{Algorithm}

\usepackage[english]{babel}

\usepackage{dsfont}
\usepackage{soul}
\usepackage{comment}
\usepackage{cases}
\usepackage{enumerate}
\usepackage[shortlabels]{enumitem}
\usepackage{tikz}
\usepackage{pgfplots}
\pgfplotsset{compat=1.15}
\usetikzlibrary{arrows}
\usetikzlibrary{patterns}
\numberwithin{equation}{section}

\usepackage{tabularx}

\newcommand{\nocontentsline}[3]{}
\let\origcontentsline\addcontentsline
\newcommand\stoptoc{\let\addcontentsline\nocontentsline}
\newcommand\resumetoc{\let\addcontentsline\origcontentsline}

\theoremstyle{plain}
\newtheorem{thm}{Theorem}[section]

\newtheorem{lem}[thm]{Lemma} 
\newtheorem{rem}[thm]{Remark} 

\newtheorem{prop}[thm]{Proposition} 
\theoremstyle{definition} 
 
\theoremstyle{remark}

\newcommand{\BB}{\mathbf{B}}

\newcommand{\Q}{\mathbf{Q}}

\newcommand{\f}{\mathbf{f}}

\newcommand{\V}{\mathbf{V}}

\makeatletter
\newcommand{\leqnomode}{\tagsleft@true}
\newcommand{\reqnomode}{\tagsleft@false}
\makeatother

\usepackage[math]{cellspace}
\title[Flux vector splittings and semi-implicit schemes for isentropic two-phase flow]{Flux vector splittings and semi-implicit finite volume schemes for an isentropic two-phase flow model}
\author{Michael Dumbser}
\author{Gabriella Puppo}
\author{Alessandra Spilimbergo}
\author{Tommaso Tenna}
\address[Michael Dumbser]{Laboratory of Applied Mathematics, DICAM, University of Trento,\newline via Mesiano 77, 38123 Trento, Italy -- \newline 
Department of Mathematics, Southern University of Science and Technology, Xueyuan Avenue 1088, 518055 Shenzhen, Guangdong, China}
\address[Gabriella Puppo]{Dipartimento di Matematica, Sapienza Universit{\`a} di Roma, P.le Aldo Moro 5, 00185 Rome, Italy}
\address[Alessandra Spilimbergo]{Dipartimento di Matematica, Sapienza Universit{\`a} di Roma, P.le Aldo Moro 5, 00185 Rome, Italy}
\address[Tommaso Tenna]{Universit{\'e} C{\^o}te d'Azur, CNRS, LJAD, Parc Valrose, F-06108 Nice, France -- Dipartimento di Matematica, Sapienza Universit{\`a} di Roma, P.le Aldo Moro 5, 00185 Rome, Italy}

\begin{document}

\maketitle

\begin{abstract}
\emergencystretch=1em
We present two new semi-implicit finite volume schemes for the one-\hspace{0pt}dimensional hyperbolic and thermodynamically compatible (HTC) isentropic two-\hspace{0pt}phase flow model of Romenski \textit{et al.}, which is written in terms of the mixture density, the volume fraction, the partial density of the gas, the mixture momentum and the relative velocity between the two phases. In gas--liquid mixtures, the sound speed of the weakly compressible liquid is much larger than the sound speed of the gas and than the flow velocity, so that the time step of fully explicit discretizations is dictated by the acoustic waves of the liquid, even when the main interest lies in the dynamics of the gas. To overcome this difficulty, the flux of the model is split at the continuous level into an explicit subsystem, which contains the dynamics of the gas, and an implicit subsystem, which contains the acoustic waves of the liquid and leads to a nonlinear system for the liquid pressure that is solved with the Newton method. In the first splitting the entire liquid phase is treated implicitly, while in the second one the convective terms of the liquid are moved to the explicit subsystem, so that only the acoustic waves of the liquid are treated implicitly and nested Picard iterations are avoided. As a consequence, the time step of both schemes is restricted only by the flow velocities and by the sound speed of the gas, and it is independent of the sound speed of the liquid.

For both schemes, we prove that states of uniform velocity and pressure across material interfaces, the so-called Abgrall condition, are preserved exactly, which is also confirmed by the numerical tests. Moreover, for well-prepared initial data, both schemes are consistent with the volume balance of the low Mach number limit of the liquid phase. A series of Riemann problems, pressure equilibrium tests and acoustic tests shows that in non-stiff regimes the semi-implicit schemes reproduce the results of a fully explicit scheme, while in the stiff regime they allow time steps that are up to three to four orders of magnitude larger, depending on the Mach number of the liquid phase, and avoid the excessive numerical dissipation of explicit schemes whose numerical viscosity is dictated by the sound speed of the liquid.
\end{abstract}

\section{Introduction}
Multi-phase flows have been widely investigated over the last decades in order to describe different types of phenomena involving interactions between different phases, such as gas--liquid flows in energy systems, particle-laden flows in environmental applications, and reactive flows in chemical engineering \cite{ishii1975, prosperetti2009, herard2019, ferrari2025}. 
Many multi-phase systems are described using diffuse interface models, in which the interface is represented as a thin transition region of finite width rather than as a sharp boundary. The seminal Baer--Nunziato model introduced in \cite{baer1986} is the cornerstone of the description of compressible two-phase flows and has inspired many alternative models, which differ in the closure relations and in the relaxation parameters, see \cite{coquel2002, dumbser2011, ferrari2024, flatten2011, gallouet2004, lund2012, saurel1999, saurel2001, saurel2018}. Among the various approaches, the hyperbolic and thermodynamically compatible (HTC) two-phase flow model of Romenski \textit{et al.} \cite{romenski2004, romenski2007, toro2010} is of particular interest because it admits a conservative formulation, which is advantageous for the theoretical analysis. The first derivation of this model is based on the theory of symmetric hyperbolic thermodynamically compatible systems \cite{godunov1995}. Since then, several studies have been devoted to the analytical and numerical investigation of this model \cite{theindumbser2022, puppo2023, thomann2023, delrio2025}, and in \cite{ferrari2025} it has been extended to an arbitrary number of phases, including viscous fluids and elastoplastic solids. In a recent work \cite{puppo2026}, the model has been derived from the multispecies Boltzmann equation, which yields explicit expressions for the coefficients and the evolution equation of the volume fraction directly from kinetic theory.

The present work focuses on the application of this framework to gas--liquid systems, which constitute one of the most relevant classes of two-phase flows in both environmental and industrial settings. The canonical example motivating this study is that of gas--liquid mixtures, in which either dispersed bubbles are immersed in and interact with a surrounding liquid, or liquid droplets are suspended in and interact with a surrounding gas. Such configurations arise, for instance, in air--water flows or in depressurizing pipe systems, where phase transitions may lead to the coexistence of liquid and gaseous states of the same medium. The two phases exhibit distinct flow properties, which may lead to significant differences in the propagation speeds of the acoustic waves and thus to different flow regimes.

For classical explicit finite volume schemes the time step is restricted by the fastest waves of both phases \cite{romenski2004, toro2010}, and schemes of Rusanov type, whose numerical viscosity is proportional to the largest signal speed, become very dissipative on the slow gas and material waves, so that very fine grids are needed to achieve acceptable accuracy. Several strategies have been proposed to overcome this difficulty, see \textit{e.g.} \cite{puppo2023}, where acoustic waves and material waves are treated separately by a semi-implicit approach. The time step is then restricted only by the material waves, which therefore remain well resolved. In the present work, instead, the gas phase is discretized entirely explicitly, and only the liquid phase, or only its acoustic waves, is treated implicitly. A different remedy, which keeps a fully explicit discretization and hence the acoustic restriction of the time step, has recently been proposed in \cite{ferrari2026} for the HTC multiphase model, where the numerical viscosity dissipates the jumps of the volume fractions only on the scale of the material velocity and the jumps of the thermodynamic state of each phase on the acoustic scale.  

The idea of treating the pressure terms implicitly, while the nonlinear convective terms are discretized explicitly, goes back to the pressure method of Casulli and Greenspan for transient compressible flows \cite{casulli1984} and to the semi-implicit finite difference schemes of Casulli for the shallow water equations \cite{casulli1990}, in which the free surface elevation is obtained from a discrete wave equation, so that the time step is no longer restricted by the celerity of the surface gravity waves. For the compressible Euler equations at low Mach number, Klein, Munz and co-workers have derived semi-implicit schemes from a multiple pressure variables asymptotic analysis, which reduce to a projection method for the incompressible Euler equations in the zero Mach number limit \cite{klein1995, munz2003}. Within this line of research, Park and Munz \cite{parkmunz2005} proposed the first conservative flow solver for all Mach numbers based on a pressure-based semi-implicit scheme on a staggered grid, in which the physical flux of the Euler equations is split at the continuous level into an advection system, discretized explicitly, and a pressure system, discretized implicitly. The same flux vector splitting was later also proposed by Toro and V{\'a}zquez-Cend{\'o}n \cite{torovazquez2012}, and it is the strategy closest to the one adopted in the present paper. When the equation of state is nonlinear, the implicit discretization of the pressure leads to nonlinear systems. If the nonlinearity is confined to the main diagonal, these systems are called mildly nonlinear, and for this class convergent Newton-type methods have been developed by Brugnano, Casulli and Zanolli \cite{brugnanocasulli2008, casullizanolli2012} and employed in conservative semi-implicit finite volume schemes for the compressible Navier--Stokes equations with general equation of state \cite{dumbsercasulli2016}. Implicit-explicit (IMEX) Runge--Kutta schemes \cite{pareschirusso2005} provide a general framework for the time integration of hyperbolic systems with stiff terms, and all speed IMEX relaxation schemes for the Euler equations have been proposed in \cite{thomann2020}.

In this paper, we consider gas--liquid mixtures in which the liquid phase is weakly compressible, \textit{i.e.} its sound speed is much larger than the sound speed of the gas and than the flow velocity, independently of whether the liquid occupies a large or a small fraction of the volume. Standard explicit methods are subject to a CFL condition dictated by the fastest waves of the system and must therefore resolve the acoustic waves of both phases \cite{toro2010}, even when the main interest lies in the dynamics of the gas. The objective of this work is therefore to design semi-implicit schemes in which the gas phase, which is not stiff, is discretized explicitly, while either the entire liquid phase or only its acoustic waves are treated implicitly. In this way the time step is no longer restricted by the sound speed of the liquid, and the schemes remain consistent with the low Mach number limit of the liquid phase.

The rest of the paper is organized as follows. In Section \ref{Section:model} we recall the isentropic two-phase flow model, the equations of state and the eigenstructure of the system, together with its formal limit for a weakly compressible liquid phase. In Section \ref{Section:flux_splitting} we introduce the two flux splittings and analyze the eigenstructure of the resulting subsystems. The two semi-implicit finite volume schemes, the preservation of the Abgrall condition and the behavior of the schemes in the low Mach number limit of the liquid are presented in Section \ref{Section:semi_implicit}, while numerical results are shown in Section \ref{Section:results}. Some concluding remarks are given in Section \ref{Section:conclusions}.

\section{The One-Dimensional HTC Isentropic Two-Phase Flow Model}
\label{Section:model}
The one-dimensional HTC isentropic two-phase flow model of Romenski \textit{et al.} \cite{romenski2004, romenski2007, toro2010} is defined by the following system

\begin{subequations}\label{system:two_phase}
\begin{eqnarray}
    \partial_t \rho + \partial_x (\rho u) &=& 0, \label{eq:mass_mix} \\
    \partial_t \alpha_1 + u \partial_x \alpha_1 &=& -\frac{1}{\tau}(p_2 - p_1), \label{eq:vol_frac} \\
    \partial_t (\alpha_1 \rho_1) + \partial_x (\alpha_1 \rho_1 u_1) &=& 0, \label{eq:mass_p1} \\
    \partial_t (\rho u) + \partial_x (\alpha_1 \rho_1 u_1^2 + \alpha_2 \rho_2 u_2^2 + \alpha_1 p_1 + \alpha_2 p_2) &=& 0, \label{eq:mom_mix} \\
    \partial_t(u_1 - u_2) + \partial_x \left( \frac{u_1^2}{2} - \frac{u_2^2}{2} + h_1 - h_2 \right) &=& -\zeta \frac{\alpha_1 \rho_1 \alpha_2 \rho_2}{\rho}(u_1 - u_2). \label{eq:rel_vel}
\end{eqnarray}
\end{subequations}

The system consists of the conservation laws for the mixture mass $\rho$ \eqref{eq:mass_mix}, for the partial mass $\alpha_1\rho_1$ \eqref{eq:mass_p1} and for the total momentum $\rho u$ \eqref{eq:mom_mix}. The evolution of the volume fraction $\alpha_1$ under the action of a pressure relaxation source term is given by \eqref{eq:vol_frac}, where $\tau$ is a relaxation parameter. Finally, the evolution of the relative velocity $u_1 - u_2$ subject to a friction source term is described by \eqref{eq:rel_vel}, with a friction coefficient $\zeta$ determined in \cite{puppo2026}. Each phase is characterized by its density, denoted by $\rho_1$ and $\rho_2$, respectively, while the mixture density $\rho$ and the mixture velocity $u$ are defined in \eqref{eq:mixture}. The volume fraction $\alpha_1 \in (0,1)$ is associated with phase one and the two volume fractions satisfy $\alpha_1 + \alpha_2 = 1$. In this paper, we use the non-conservative form \eqref{eq:vol_frac} of the volume fraction equation instead of the conservative form $\partial_t(\rho\alpha_1)+\partial_x(\rho\alpha_1 u)=-\frac{\rho}{\tau}(p_2-p_1)$. The two formulations are equivalent only for smooth solutions, and in our numerical experiments the conservative form led to considerably worse results than the non-conservative one, very likely because the volume fraction is then recovered as $\alpha_1=(\rho\alpha_1)/\rho$, so that the numerical diffusion acting on the mixture density $\rho$ enters the volume fraction. The non-conservative form, in contrast, allows the preservation of pressure--velocity equilibria across material interfaces, see \cite{abgrall1996, saurel1999} and Section \ref{Section:semi_implicit}. 

\subsection{Mixture Definitions and Equations of State}
The mixture density $\rho$ and velocity $u$ are defined as \cite{toro2010}
\begin{equation}
\label{eq:mixture}
    \rho = \alpha_1 \rho_1 + \alpha_2 \rho_2, \quad u = \frac{\alpha_1 \rho_1 u_1 + \alpha_2 \rho_2 u_2}{\rho},
\end{equation}
and we denote by $\chi_k = \alpha_k \rho_k / \rho$, $k=1,2$, the mass fractions of the two phases. The phase enthalpies are $h_i = e_i + \frac{p_i}{\rho_i}$ \cite{puppo2023}. To close the system, let us introduce an equation of state for each phase. We use the isentropic ideal gas law for phase 1 and the stiffened gas law for phase 2 \cite{harlow1970, saurel1999},
namely 
\begin{eqnarray}
    e_1(\rho_1) &=& \frac{C_1^2}{\gamma_1 (\gamma_1-1)} \left( \frac{\rho_1}{\rho_{01}} \right)^{\gamma_1-1} \\
    e_2(\rho_2) &=& \frac{C_2^2}{\gamma_2 (\gamma_2-1)} \left( \frac{\rho_2}{\rho_{02}} \right)^{\gamma_2-1} + \frac{p_\infty}{\rho_2}
\end{eqnarray}
from which the pressures are computed via $p_i = \rho_i^2\, \partial e_i / \partial \rho_i$ as
\begin{eqnarray}
    p_1(\rho_1) &=& \frac{\rho_{01} C_1^2}{\gamma_1} \left( \frac{\rho_1}{\rho_{01}} \right)^{\gamma_1} \label{eq:eos_p1} \\
    p_2(\rho_2) &=& \frac{\rho_{02} C_2^2}{\gamma_2} \left( \frac{\rho_2}{\rho_{02}} \right)^{\gamma_2} - p_\infty, \label{eq:eos_p2}
\end{eqnarray}
where $\rho_{0i}$ is the reference density, $\gamma_i$ the adiabatic exponent and $C_i$ the sound speed of phase $i$ at normal atmospheric conditions, and $p_\infty = (\rho_{02} C_2^2 - \gamma_2 p_0)/\gamma_2$. Since $p_2 + p_\infty > 0$ for every $\rho_2 > 0$, the liquid pressure may also take negative values, and \eqref{eq:eos_p2} can be inverted to $\rho_2 = \rho_2(p_2)$ for all $p_2 > -p_\infty$.
In this work, the reference pressure of the stiffened gas is chosen equal to the pressure of the ideal gas at $\rho_1 = \rho_{01}$, namely
\begin{equation}
    p_0 = \frac{\rho_{01} C_1^2}{\gamma_1}.
\end{equation}
Let us now analyze the eigenstructure of the system, which can be written in compact form 
\begin{equation}
\label{eq:compact_form_1}
    \partial_t \Q + \partial_x \f(\Q) + \BB(\Q) \partial_x \Q =\mathbf{S}(\Q)
\end{equation}
by introducing the state vector 
\begin{equation}
\label{eq:state_vector}
    \Q = (\rho, \alpha_1, \alpha_1 \rho_1, \rho u, u_1-u_2)^T,
\end{equation}
the physical conservative flux $\f$ of the system
\begin{equation}
    \renewcommand{\arraystretch}{1.3}
    \f = \begin{pmatrix}
    \rho u\\ 
    0\\
    \alpha_1 \rho_1 u_1 \\
    \alpha_1 \rho_1 u_1^2 + \alpha_2 \rho_2 u_2^2 + \alpha_1 p_1 + \alpha_2 p_2\\
    \dfrac{u_1^2}{2} - \dfrac{u_2^2}{2} + h_1 - h_2
    \end{pmatrix},
\end{equation}
$\BB(\Q)=u\,\mathbf{e}_2\mathbf{e}_2^T$ with $\mathbf{e}_2=(0,1,0,0,0)^T$, and the source term
\begin{equation}
    \renewcommand{\arraystretch}{1.3}
    \mathbf{S} = \begin{pmatrix}
    0\\
    -\dfrac{1}{\tau}(p_2-p_1)\\
    0 \\
    0 \\
    - \zeta \dfrac{\alpha_1 \rho_1 \alpha_2 \rho_2}{\rho} (u_1-u_2)
    \end{pmatrix}.
\end{equation}
The model is hyperbolic with the following eigenvalues
\begin{equation}
\label{eq:eigenvalues_complete_system}
    \lambda^u = u, \qquad \lambda_1^{\pm} = u_1 \pm c_1, \qquad  \lambda_2^{\pm} = u_2 \pm c_2,
\end{equation}
which are the eigenvalues of the system matrix $\mathbf{A}(\Q) = \mathbf{J}(\Q) + \BB(\Q)$ with the Jacobian matrix $\mathbf{J} = \partial \f / \partial \Q$, and where $c_i^2 = \tfrac{\partial p_i}{\partial \rho_i}$. Therefore, the dynamics of the system is governed by two acoustic waves per phase $\lambda_{1,2}^{\pm}$ and the mixture velocity $\lambda^u$, referred to as \textit{material velocity}. This means that for large sound speeds $c_{1,2}$, the acoustic waves $\lambda_{1,2}^\pm$ travel much faster than the material wave $\lambda^u$, leading to multiscale dynamics. In particular, we can introduce a reference quantity called \textit{Mach number}, given by the ratio between the phase velocities and the sound speeds, namely
\begin{equation}
    M_i = \frac{u_i}{c_i}, \qquad i=1,2.
\end{equation}
For low Mach numbers, the acoustic waves propagate significantly faster than the material wave. In the following, we are only interested in two-phase flows in which one of the two phases is characterized by a low Mach number. 

\subsection{Asymptotic Analysis in the Limit $M_2 \ll M_1=1$}
Let the gas phase (phase 1) be compressible with $\rho_1 = \mathcal{O}(1)$, $c_1 = \mathcal{O}(1)$ and $u_1 = \mathcal{O}(1)$. 
We expand the variables of the weakly compressible liquid phase (phase 2) with respect to the Mach number $M_2$ as follows
\begin{equation}
\label{eq:R2_U2_exp}
    \rho_2 = \rho_2^{(0)} + M_2 \rho_2^{(1)} + M_2^2 \rho_2^{(2)} + \mathcal{O}(M_2^3), \quad u_2 = u_2^{(0)} + M_2 u_2^{(1)} + \mathcal{O}(M_2^2).
\end{equation}
Using the EOS \eqref{eq:eos_p2}, we obtain the Mach number expansion of the pressure from the density expansion \eqref{eq:R2_U2_exp}. In particular, with $p_2^{(0)} = p_2(\rho_2^{(0)})$ and $(c_2^{(0)})^2 = \frac{\partial p_2}{\partial \rho_2}(\rho_2^{(0)})$, we get
\begin{equation}
p_2 = p_2^{(0)} + (c_2^{(0)})^2\rho_2^{(1)} M_2 + \left((c^{(0)}_2)^2\rho_2^{(2)} + \frac{1}{2}(\gamma_2-1)\frac{(c^{(0)}_2)^2}{\rho_2^{(0)}} (\rho_2^{(1)})^2\right) M_2^2 + \mathcal{O}(M_2^3). 
\end{equation}
Furthermore, let us assume the following Mach number expansion of the volume fraction
\begin{equation}
\label{eq:M_exp_alpha}
    \alpha = \alpha^{(0)} + \mathcal{O}(M_2),
\end{equation}
where here and in the following $\alpha = \alpha_1$. As in \cite{puppo2023}, we write \eqref{system:two_phase} in non-dimensional form. Both phases share a common reference velocity $u_r = x_r/t_r$ and a reference density $\rho_r$. The reference pressure is $\rho_r u_r^2$ for the gas and $\rho_r c_{2,r}^2 = \rho_r u_r^2/M_2^2$ for the liquid, where $c_{2,r}$ is a reference sound speed of the liquid and $M_2 = u_r/c_{2,r}$. The reference values for $\tau$ and $\zeta$ are $\rho_r u_r^2 t_r$ and $1/(\rho_r t_r)$. The ratios of the reference phase densities to $\rho_r$ are treated as constants of order one and omitted. The scaled system is then \eqref{system:two_phase} with $p_2$ and $h_2$ replaced by $p_2/M_2^2$ and $h_2/M_2^2$, and all quantities in the remainder of this section are non-dimensional. The following result holds
\begin{lem}[\cite{puppo2023}, Well-prepared data for compressible/weakly compressible flow]
    \label{lem:well_prepared}
    Consider a compressible gas phase (phase 1) with $M_1 = 1$ and a weakly compressible liquid phase (phase 2) with $M_2 = M \ll 1$, and assume that $\alpha$, $\rho_2$ and $u_2$ admit the Mach number expansions \eqref{eq:M_exp_alpha} and \eqref{eq:R2_U2_exp}. We call well-prepared the initial data in the set
    \begin{equation}
    	\label{eq:wp_1}
    	\Omega_1^{wp} = \left\lbrace \Q \in \mathbb{R}^{5}: \rho_2^{(0)} \text{ const.},~ p_2^{(0)} = 0, ~\rho_2^{(1)} = 0, ~\partial_t \alpha^{(0)} = 0, ~\partial_x u_2^{(0)} = \frac{u_2^{(0)}}{\alpha_2^{(0)}} \partial_x \alpha^{(0)} \right\rbrace.
	\end{equation}
	Then, formally, for $\Q \in \Omega_1^{wp}$ and $M \to 0$, the limit equations are given by
	\begin{subequations}\label{system:two_phase_asymptotic_Mach}
\begin{eqnarray}
    & & u^{(0)} \partial_x \alpha^{(0)} = - \frac{1}{\tau}(p_2^{(2)} - p_1), \label{eq:Limit1_alpha} \\
    & & \partial_t u_2^{(0)} + u_2^{(0)}\partial_x u_2^{(0)} + \frac{\partial_x p_2^{(2)}}{\rho_2^{(0)}} + \frac{p_1 - p_2^{(2)}}{\rho^{(0)}}\partial_x \alpha^{(0)} = \zeta \rho^{(0)} (\chi_1^{(0)})^2\chi_2^{(0)} (u_1 - u_2^{(0)}), \label{eq:Limit1_u2} \\
    & & \partial_t (\alpha^{(0)} \rho_1) + \partial_x (\alpha^{(0)}\rho_1 u_1) = 0, \label{eq:Limit1_ar1} \\
    & & \partial_t u_1 + u_1 \partial_x u_1 + \frac{\partial_x p_1}{\rho_1} + \frac{p_1 - p_2^{(2)}}{\rho^{(0)}}\partial_x \alpha^{(0)} = -\zeta \rho^{(0)} \chi_1^{(0)}(\chi_2^{(0)})^2 (u_1 - u_2^{(0)}). \label{eq:Limit1_u1}
\end{eqnarray}
\end{subequations}  
	If, in addition, $\partial_x \alpha^{(0)} = 0$, \textit{i.e.} $\alpha^{(0)}$ is constant, the limit equations reduce to 
	\begin{equation}
		\label{eq:Limit1_cons}
    	\begin{split}
    	\partial_t u_2^{(0)} + \frac{\partial_x p_1}{\rho_2^{(0)}} &= \zeta \rho^{(0)} (\chi_1^{(0)})^2 \chi_2^{(0)} (u_1 - u_2^{(0)}), \quad \partial_x u_2^{(0)} = 0,\\
    	\partial_t \rho_1 + \partial_x (\rho_1 u_1) & = 0, \\
    	\partial_t (\rho_1 u_1) + \partial_x (\rho_1 u_1^2 + p_1) &= - \zeta \rho^{(0)} \rho_1 \chi_1^{(0)}(\chi_2^{(0)})^2 (u_1 - u_2^{(0)}).
    	\end{split}
	\end{equation}
\end{lem}
    	
The proof proceeds as in \cite[Section 2.3]{puppo2023}, taking into account that $\tau$ and $\zeta$ are defined differently there. In the limit $M \to 0$, phase two is therefore governed by the incompressible Euler equations with friction and phase one by the compressible isentropic Euler equations with friction. The two phases remain coupled through the friction terms and through the gas pressure $p_1$, which is the only pressure left in the limit and also drives the liquid. Note that \eqref{system:two_phase_asymptotic_Mach} and \eqref{eq:Limit1_cons} refer to the model with pressure relaxation: for constant $\alpha^{(0)}$, it is \eqref{eq:Limit1_alpha} with a finite relaxation time $\tau$ that yields $p_2^{(2)} = p_1$. The schemes of Section~\ref{Section:semi_implicit} solve only the homogeneous part of \eqref{system:two_phase}, and the consistency argument of Section~\ref{sec:low_mach_consistency} uses only the liquid volume balance, which holds also without relaxation.

\section{The Flux-Splitting Approach}
\label{Section:flux_splitting}

\begin{figure}[htbp]
    \centering
    \includegraphics[width=\textwidth]{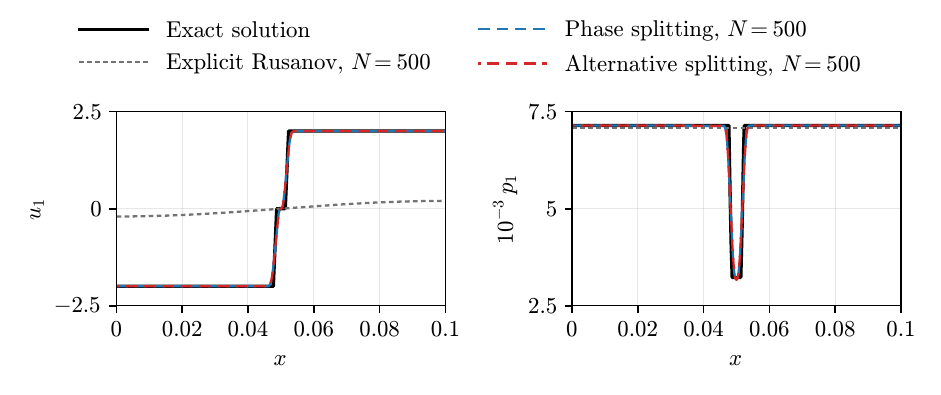}
    \caption{Numerical approximation of a gas--liquid mixture with a liquid volume fraction of $\alpha_2=0.03$: gas velocity $u_1$ (left) and gas pressure $p_1$ (right), compared with the exact solution \cite{theindumbser2022}, see Section~\ref{test_1}. The comparison shows that the explicit scheme is significantly more diffusive, since its numerical viscosity is dictated by the sound speed of the liquid.}
    \label{fig:test_2000}
\end{figure}

Since the goal is to construct a semi-implicit scheme in which $\Delta t$ is determined by a CFL condition based on the eigenvalues of the \textit{less stiff} phase, \textit{i.e.} the one with the larger Mach number, in this section we propose a splitting at the continuous level. The proposed strategy is based on a flux splitting, in which the dynamics of the two phases are separated. The motivation for such a splitting is illustrated in Figure \ref{fig:test_2000}. A standard explicit scheme for the full system requires a time step $\Delta t \approx 10^{-9}$ for stability, even when the liquid phase constitutes only $3\%$ of the total volume, whereas a semi-implicit scheme based on the splitting approach detailed in this section allows a much larger time step, $\Delta t \approx 10^{-5}$, dictated only by the gasdynamics. Moreover, since the numerical viscosity of the explicit Rusanov scheme is dictated by the sound speed of the liquid, the explicit solution is strongly smeared, as shown in Figure \ref{fig:test_2000}. 
For this test, the computational domain is $[0,L]$ with
$L=0.1$, discretized with $N_x=500$ cells and $\nu_{\text{CFL}}=0.6$, with transmissive boundary conditions and without source terms.
The parameters of the equations of state are $\gamma_1=1.4$, $\rho_{01}=1$, $C_1=1$
for the gas and $\gamma_2=2.8$, $\rho_{02}=10^3$, $C_2=5\cdot10^4$ for the liquid. The initial data are $\rho_1=719.68$, $\rho_2=1320$, $\alpha_1=0.97$ in the whole
domain and $u_1=u_2=-2$ for $x<L/2$, $u_1=u_2=2$ for $x\geq L/2$, while the final time is
$T=4\cdot10^{-4}$. The semi-implicit solution is computed with the two splitting schemes presented in this paper, which require only 20 time steps with
$\Delta t\approx2.1\cdot10^{-5}$, while the explicit Rusanov scheme requires about $2.1\cdot10^{5}$ time steps with $\Delta t\approx1.9\cdot10^{-9}$. These results serve as a motivation of the work proposed in this paper.

\subsection{Phase-Splitting Scheme}
\label{Section:phase_splitting}

Let us recall the compact form \eqref{eq:compact_form_1}
\begin{equation*}
    \partial_t \Q + \partial_x \f(\Q) + \BB(\Q) \partial_x \Q =\mathbf{S}(\Q)
\end{equation*}
for the state vector $\Q$. The strategy consists in splitting the flux $\f$ into two conservative fluxes $\f_1$ and $\f_2$ such that
\begin{equation}
    \f(\Q) = \f_1(\Q) + \f_2(\Q),
\end{equation}
in order to rewrite \eqref{eq:compact_form_1} as
\begin{equation}
    \partial_t \Q + \partial_x \f_1(\Q) + \partial_x \f_2(\Q) + \BB(\Q) \partial_x \Q = \mathbf{S}(\Q),
\end{equation}
where
\begin{equation}
    \renewcommand{\arraystretch}{1.3}
\label{eq:fluxes_phase}
    \f_1 = \begin{pmatrix}
    \alpha_1 \rho_1 u_1\\ 
    0\\
    \alpha_1 \rho_1 u_1 \\
    \alpha_1 \rho_1 u_1^2 + \alpha_1 p_1\\
    \dfrac{u_1^2}{2} + h_1
    \end{pmatrix}, \qquad     \f_2 = \begin{pmatrix}
    \alpha_2 \rho_2 u_2\\ 
    0\\
    0\\
    \alpha_2 \rho_2 u_2^2 + \alpha_2 p_2\\
    -\dfrac{u_2^2}{2} - h_2
    \end{pmatrix}.
\end{equation}
From this flux decomposition, we propose the following phase-splitting into two subsystems
\begin{subequations}
\begin{eqnarray}
    \partial_t \Q + \partial_x \f_1(\Q) + \BB(\Q) \partial_x \Q &=& 0, \label{eq:compact_subystem1} \\
    \partial_t \Q + \partial_x \f_2(\Q) &=& \mathbf{S}(\Q). \label{eq:compact_subystem2}
\end{eqnarray}
\end{subequations}
In the following, we refer to subsystem \eqref{eq:compact_subystem1} as the gas-system and to subsystem \eqref{eq:compact_subystem2} as the liquid-system. For $i = 1, 2$, we denote by $\mathbf{J}_i = \partial \f_i / \partial \Q$ the Jacobian matrices of the two fluxes, by $\BB_1 = \BB(\Q)$ and $\BB_2 = 0$ the matrices of the non-conservative products, and by $\mathbf{A}_i = \mathbf{J}_i + \BB_i$ the system matrices of the two subsystems.

\subsubsection{The gas subsystem}
The gas-system is given by
\begin{subequations}\label{eq:gas_system}
\begin{eqnarray}
    \partial_t \rho + \partial_x (\alpha_1 \rho_1 u_1) &=& 0, \label{eq:gas-system:1} \\
    \partial_t \alpha_1 + u \partial_x \alpha_1 &=& 0, \label{eq:gas-system:2} \\
    \partial_t (\alpha_1 \rho_1) + \partial_x (\alpha_1 \rho_1 u_1) &=& 0, \label{eq:gas-system:3} \\
    \partial_t (\rho u) + \partial_x (\alpha_1 \rho_1 u_1^2 + \alpha_1 p_1) &=& 0, \label{eq:gas-system:4} \\
    \partial_t (u_1-u_2) + \partial_x \left(\frac{u_1^2}{2} + h_1\right) &=& 0. \label{eq:gas-system:5}
\end{eqnarray}
\end{subequations}
\begin{prop}
\label{prop:gas_system}
System \eqref{eq:compact_subystem1} has the eigenvalues
\begin{equation}
    \lambda_{1,1}=u_1+c_1, \quad \lambda_{1,2}=u_1-c_1, \quad \lambda_{1,3}=u, \quad \lambda_{1,4}=0, \quad \lambda_{1,5}=0,
\end{equation}
and is strongly hyperbolic, except for $u = 0$ or $(u_1-u)^2 = c_1^2$, where it is only weakly hyperbolic.
\end{prop}
\begin{proof}
    Let us compute the eigenvalues in terms of the vector of primitive variables $\V = \left[\alpha_1, \rho_1, u_1, \rho_2, u_2 \right]^T$.
    Equation \eqref{eq:gas-system:2} is already given for $\alpha_1$.

    Equation \eqref{eq:gas-system:3} may be recast as
    \begin{equation}
        \alpha_1 \partial_t \rho_1 + \rho_1 \partial_t \alpha_1 + \alpha_1 \rho_1 \partial_x u_1 + u_1 \alpha_1 \partial_x \rho_1 + u_1 \rho_1 \partial_x \alpha_1 = 0,
    \end{equation}
    and using \eqref{eq:gas-system:2} we obtain
    \begin{equation}
        \partial_t \rho_1 + \rho_1 \partial_x u_1 + u_1 \partial_x \rho_1 + \frac{\rho_1}{\alpha_1}(u_1-u) \partial_x \alpha_1 = 0.
    \end{equation}
    Let us now construct the equation for $u_1$. Starting from Equation \eqref{eq:gas-system:4}, we can write
    \begin{equation}
        \alpha_1\rho_1\partial_t(u_1) + \alpha_2\rho_2\partial_t(u_2) + (\alpha_1\rho_1u_1)\partial_x(u_1) + \partial_x(\alpha_1 p_1) = 0,
    \end{equation}
    which, combined with \eqref{eq:gas-system:5}, leads to
    \begin{equation}
        \partial_t(u_1) + \frac{p_1}{\rho} \partial_x \alpha_1 + \frac{c_1^2}{\rho_1} \partial_x (\rho_1) + u_1 \partial_x u_1 = 0.
    \end{equation}
    Equation \eqref{eq:gas-system:1} can be rewritten as
    \begin{equation}
        \partial_t (\alpha_1 \rho_1 + \alpha_2 \rho_2) + \partial_x (\alpha_1 \rho_1 u_1) = 0,
    \end{equation}
    which, using \eqref{eq:gas-system:3}, leads to
    \begin{equation}
        \partial_t (\alpha_2 \rho_2) = 0.
    \end{equation}
    This is equivalent to 
    \begin{equation}
        \alpha_2 \partial_t \rho_2 - \rho_2 \partial_t \alpha_1 = 0 \implies \partial_t \rho_2 + \frac{\rho_2 u}{\alpha_2} \partial_x \alpha_1 = 0.
    \end{equation}
     Finally, the equation for $u_2$ becomes
     \begin{equation}
         \partial_t(u_2) + \frac{p_1}{\rho} \partial_x \alpha_1 = 0.
     \end{equation}
     Therefore, in terms of the primitive variables the gas-system reads $\partial_t \V + \mathbf{C}_1  \, \partial_x \V = 0$, where
     \begin{equation}
    \renewcommand{\arraystretch}{1.3}
    \label{eq:C_gas}
         \mathbf{C}_1 = \frac{\partial \V}{\partial \Q} \, \mathbf{A}_1 \, \frac{\partial \Q}{\partial \V} = \begin{pmatrix}
         u & 0 & 0 & 0 & 0\\
         \dfrac{\alpha_2 \rho_2 \rho_1}{\alpha_1 \rho} (u_1-u_2) & u_1 & \rho_1 & 0 & 0\\
         \dfrac{p_1}{\rho} & \dfrac{c_1^2}{\rho_1} & u_1 & 0 & 0\\
         \dfrac{\rho_2 u}{\alpha_2} & 0 & 0 & 0 & 0\\
         \dfrac{p_1}{\rho} & 0 & 0 & 0 & 0
         \end{pmatrix},
     \end{equation}
     from which the eigenvalues can be easily computed.
    Moreover, by solving $(\mathbf{C}_1 - \lambda \mathbf{I})\mathbf{r} = 0$, we obtain the following complete set of linearly independent eigenvectors.

    For the acoustic waves $\lambda_{1,1}, \lambda_{1,2} = u_1 \pm c_1$ we get    
    \begin{equation}
        \mathbf{r}_1 = \begin{pmatrix} 0 \\ \rho_1 \\ c_1 \\ 0 \\ 0 \end{pmatrix}, \quad 
        \mathbf{r}_2 = \begin{pmatrix} 0 \\ \rho_1 \\ -c_1 \\ 0 \\ 0 \end{pmatrix}.
    \end{equation}
    
    For the interfacial wave $\lambda_{1,3} = u$, defining $D = (u_1-u)^2 - c_1^2$ and assuming $u \neq 0$ and $D \neq 0$, the corresponding eigenvector is:
    \begin{equation}
    \renewcommand{\arraystretch}{1.3}
        \mathbf{r}_3 = \begin{pmatrix} 
        1 \\ 
        \dfrac{\rho_1 C_{31} - C_{21}(u_1-u)}{D} \\ 
        \dfrac{C_{21} (c_1^2/\rho_1) - C_{31}(u_1-u)}{D} \\ 
        \dfrac{\rho_2}{\alpha_2} \\ 
        \dfrac{p_1}{u \rho} 
        \end{pmatrix}
    \end{equation}
    where $C_{21}$ and $C_{31}$ denote the entries of $\mathbf{C}_1$ in \eqref{eq:C_gas}.
    
    For the stationary waves $\lambda_{1,4} = \lambda_{1,5} = 0$, the geometric multiplicity is 2, which yields two linearly independent eigenvectors, namely the fourth and fifth standard basis vectors:
    \begin{equation}
        \mathbf{r}_4 = \begin{pmatrix} 0 \\ 0 \\ 0 \\ 1 \\ 0 \end{pmatrix}, \quad 
        \mathbf{r}_5 = \begin{pmatrix} 0 \\ 0 \\ 0 \\ 0 \\ 1 \end{pmatrix}
    \end{equation}
    
    Since this system has five real eigenvalues and a complete set of five linearly independent eigenvectors, the system is strongly hyperbolic, except at the singular points $u = 0$ and $(u_1 - u)^2 = c_1^2$, where it is only weakly hyperbolic.
\end{proof}

\begin{rem}
Three of the five eigenvalues of system \eqref{eq:compact_subystem1} coincide with those of the complete system \eqref{eq:eigenvalues_complete_system}, namely
\begin{equation}
    \lambda_{1,1} = \lambda_1^+, \quad \lambda_{1,2} = \lambda_1^-, \quad \lambda_{1,3} = \lambda^u.
\end{equation}

\end{rem}

\subsubsection{The liquid subsystem}
The liquid-system is given by
\begin{subequations}\label{eq:liquid-system}
\begin{eqnarray}
    \partial_t \rho + \partial_x (\alpha_2 \rho_2 u_2) &=& 0, \\
    \partial_t \alpha_1 &=& -\frac{1}{\tau}(p_2-p_1), \\
    \partial_t (\alpha_1 \rho_1) &=& 0, \\
    \partial_t (\rho u) + \partial_x (\alpha_2 \rho_2 u_2^2 + \alpha_2 p_2) &=& 0, \\
    \partial_t (u_1-u_2) - \partial_x \left(\frac{u_2^2}{2} + h_2\right) &=& -\zeta \frac{\alpha_1 \rho_1 \alpha_2 \rho_2}{\rho}(u_1 - u_2).
\end{eqnarray}
\end{subequations}

The liquid-system \eqref{eq:liquid-system} has the eigenvalue $0$ with algebraic multiplicity three and the acoustic eigenvalues $\lambda = u_2 \pm c_2$. Since $\partial_x(\alpha_2 p_2) = \alpha_2 \partial_x p_2 - p_2 \partial_x \alpha_1$, the equations for both phase velocities contain the term $-\frac{p_2}{\rho}\partial_x \alpha_1$, so that for $p_2 \neq 0$ the eigenvalue $0$ has only two linearly independent eigenvectors and the liquid-system is only weakly hyperbolic. The same holds for the explicit convective subsystem of the flux vector splitting of  Toro and V{\'a}zquez-Cend{\'o}n \cite{torovazquez2012}, whose double eigenvalue $u$ has only one eigenvector.

\subsection{Alternative Splitting}
\label{Section:alternative_splitting}
The phase-splitting approach simplifies the dynamics by treating the gas phase, which is much slower than the liquid phase, explicitly. However, treating the entire liquid phase implicitly requires the solution of the nonlinear discrete wave equation \eqref{eq:wave_phase}, which also contains the convective terms of the liquid and therefore nested Picard iterations, see Section \ref{Section:numerical_phase_splitting}.

To this end, following \cite{puppo2023}, we propose an alternative splitting, in which the gas phase is still treated fully explicitly, but also the convective part of the liquid phase is treated explicitly. In other words, in the liquid phase we separate the material waves from the acoustic waves, and only the latter are treated implicitly.

Let us introduce the following splitting of $\f$ into two conservative fluxes $\hat{\f}_1$ and $\hat{\f}_2$ such that
\begin{equation}
    \f(\Q) = \hat{\f}_1(\Q) + \hat{\f}_2(\Q),
\end{equation}
where now
\begin{equation}
    \renewcommand{\arraystretch}{1.3}
     \hat{\f}_1 = \begin{pmatrix}
    \alpha_1 \rho_1 u_1\\ 
    0\\
    \alpha_1 \rho_1 u_1 \\
    \alpha_1 \rho_1 u_1^2 + \alpha_1 p_1 + \alpha_2 \rho_2 u_2^2\\
    \dfrac{u_1^2}{2} + h_1 -  \dfrac{u_2^2}{2}
    \end{pmatrix}, \qquad     
    \hat{\f}_2 = \begin{pmatrix}
    \alpha_2 \rho_2 u_2\\ 
    0\\
    0\\
    \alpha_2 p_2\\
    -h_2
    \end{pmatrix}.
\end{equation}
From this flux decomposition, we obtain the following splitting into two subsystems:
\begin{subequations}
\begin{eqnarray}
    \partial_t \Q + \partial_x \hat{\f}_1(\Q) + \BB(\Q) \partial_x \Q &=& 0, \label{eq:compact_subystem1_alternative} \\
    \partial_t \Q + \partial_x \hat{\f}_2(\Q) &=& \mathbf{S}(\Q). \label{eq:compact_subystem2_alternative}
\end{eqnarray}
\end{subequations}
In the following, we refer to \eqref{eq:compact_subystem1_alternative} as the first subsystem and to \eqref{eq:compact_subystem2_alternative} as the second subsystem. Analogously, we denote by $\hat{\mathbf{J}}_i = \partial \hat{\f}_i / \partial \Q$ the Jacobian matrices of the two fluxes and by $\hat{\mathbf{A}}_i = \hat{\mathbf{J}}_i + \BB_i$ the system matrices of the first subsystem and of the second subsystem, with $\BB_1 = \BB(\Q)$ and $\BB_2 = 0$ as before.

\subsubsection{The first subsystem}
The first subsystem is given by
\begin{subequations}\label{eq:1-system}
\begin{eqnarray}
    \partial_t \rho + \partial_x (\alpha_1 \rho_1 u_1) &=& 0, \label{eq:1-system:1} \\
    \partial_t \alpha_1 + u \partial_x \alpha_1 &=& 0, \label{eq:1-system:2} \\
    \partial_t (\alpha_1 \rho_1) + \partial_x (\alpha_1 \rho_1 u_1) &=& 0, \label{eq:1-system:3} \\
    \partial_t (\rho u) + \partial_x (\alpha_1 \rho_1 u_1^2 + \alpha_1 p_1 + \alpha_2 \rho_2 u_2^2) &=& 0, \label{eq:1-system:4} \\
    \partial_t (u_1-u_2) + \partial_x \left(\frac{u_1^2}{2} + h_1 - \frac{u_2^2}{2}\right) &=& 0. \label{eq:1-system:5}
\end{eqnarray}
\end{subequations}

\begin{prop}
System \eqref{eq:compact_subystem1_alternative} has the eigenvalues 
\begin{equation}
     \lambda_{1,1}=u_1+c_1, \quad \lambda_{1,2}= u_1-c_1, \quad \lambda_{1,3}=u, \quad \lambda_{1,4}= u_2 \left(1+\frac{\alpha_2 \rho_2}{\rho} \right), \quad \lambda_{1,5}=0,
\end{equation}
and is strongly hyperbolic, except for $u = \lambda_{1,4}$, $(u_1-u)^2 = c_1^2$, $u_1^2 = c_1^2$ or $(u_1-\lambda_{1,4})^2 = c_1^2$, where for $u_2 \neq 0$ it is in general only weakly hyperbolic. For $u_2 = 0$, the matrix $\hat{\mathbf{C}}_1$ coincides with the matrix $\mathbf{C}_1$ of the gas subsystem, and Proposition~\ref{prop:gas_system} applies.
\end{prop}

\begin{proof}
    Let us compute the eigenvalues with respect to the vector of primitive variables $\V = \left[\alpha_1, \rho_1, u_1, \rho_2, u_2 \right]^T$.
    Equation \eqref{eq:1-system:2} is already given for $\alpha_1$.

    Equation \eqref{eq:1-system:3} may be recast as
    \begin{equation}
        \alpha_1 \partial_t \rho_1 + \rho_1 \partial_t \alpha_1 + \alpha_1 \rho_1 \partial_x u_1 + u_1 \alpha_1 \partial_x \rho_1 + u_1 \rho_1 \partial_x \alpha_1 = 0,
    \end{equation}
    and using \eqref{eq:1-system:2} we obtain
    \begin{equation}
        \partial_t \rho_1 + \rho_1 \partial_x u_1 + u_1 \partial_x \rho_1 + \frac{\rho_1}{\alpha_1}(u_1-u) \partial_x \alpha_1 = 0.
    \end{equation}
    Equation \eqref{eq:1-system:1} can be rewritten as
    \begin{equation}
        \partial_t (\alpha_1 \rho_1 + \alpha_2 \rho_2) + \partial_x (\alpha_1 \rho_1 u_1) = 0,
    \end{equation}
    which, using \eqref{eq:1-system:3}, yields
    \begin{equation}
        \partial_t (\alpha_2 \rho_2) = 0.
    \end{equation}
    This is equivalent to 
    \begin{equation}
        \alpha_2 \partial_t \rho_2 - \rho_2 \partial_t \alpha_1 = 0 \implies \partial_t \rho_2 + \frac{\rho_2 u}{\alpha_2} \partial_x \alpha_1 = 0.
    \end{equation}
     Let us construct the equations for $u_1$ and $u_2$. Expanding the convective derivative in \eqref{eq:1-system:4} gives
     \begin{equation}
         \alpha_1 \rho_1 (\partial_t u_1 + u_1 \partial_x u_1) + \alpha_2 \rho_2 (\partial_t u_2 + u_2 \partial_x u_2) + u_2 \partial_x (\alpha_2 \rho_2 u_2) + \partial_x (\alpha_1 p_1) = 0.
     \end{equation}
     Now we can rewrite the relative velocity equation \eqref{eq:1-system:5} recalling that $\partial_x h_k = \frac{1}{\rho_k} \partial_x p_k$, thus
     \begin{equation}
         \partial_t u_1 - \partial_t u_2 + u_1 \partial_x u_1 - u_2 \partial_x u_2 + \frac{1}{\rho_1} \partial_x p_1 = 0.
     \end{equation}
     For the liquid phase, this gives
     \begin{equation}
        \partial_t u_2 + u_2 \partial_x u_2 = (\partial_t u_1 + u_1 \partial_x u_1) + \frac{1}{\rho_1} \partial_x p_1.
     \end{equation}
     By substituting this expression into the equation for the total momentum, we have
     \begin{multline}
         \alpha_1 \rho_1 (\partial_t u_1 + u_1 \partial_x u_1) + \alpha_2 \rho_2 \left[ (\partial_t u_1 + u_1 \partial_x u_1) + \frac{1}{\rho_1} \partial_x p_1 \right] \\ + u_2^2 \alpha_2 \partial_x \rho_2 + u_2^2 \rho_2 \partial_x \alpha_2 + u_2 \alpha_2 \rho_2 \partial_x u_2+ \alpha_1 \partial_x p_1 + p_1 \partial_x \alpha_1 = 0,
     \end{multline}
     which is equivalent to
     \begin{multline}
        (\alpha_1 \rho_1 + \alpha_2 \rho_2)(\partial_t u_1 + u_1 \partial_x u_1) + \alpha_2 \frac{\rho_2}{\rho_1} \partial_x p_1 \\ + u_2^2 \alpha_2 \partial_x \rho_2 + u_2^2 \rho_2 \partial_x \alpha_2 + u_2 \alpha_2 \rho_2 \partial_x u_2 + \alpha_1 \partial_x p_1 + p_1 \partial_x \alpha_1 = 0.
     \end{multline}
     We observe that $\alpha_1 \rho_1 + \alpha_2 \rho_2 = \rho$ and that the coefficient of $\partial_x p_1$ is $\alpha_1 + \alpha_2 \frac{\rho_2}{\rho_1} = \frac{\alpha_1 \rho_1 + \alpha_2 \rho_2}{\rho_1} = \frac{\rho}{\rho_1}$.

     Substituting and dividing by $\rho$, we obtain the primitive equation for the gas velocity
     \begin{equation}
         \partial_t u_1 + u_1 \partial_x u_1 + \frac{c_1^2}{\rho_1} \partial_x \rho_1 + \frac{p_1}{\rho} \partial_x \alpha_1 + \frac{1}{\rho}\left(u_2^2 \alpha_2 \partial_x \rho_2 + u_2^2 \rho_2 \partial_x \alpha_2 + u_2 \alpha_2 \rho_2 \partial_x u_2\right) = 0,
     \end{equation}
     and, substituting back, the equation for the liquid velocity
     \begin{equation}
        \partial_t u_2 + u_2 \partial_x u_2 + \frac{p_1}{\rho} \partial_x \alpha_1 + \frac{1}{\rho}\left(u_2^2 \alpha_2 \partial_x \rho_2 + u_2^2 \rho_2 \partial_x \alpha_2 + u_2 \alpha_2 \rho_2 \partial_x u_2\right) = 0.
     \end{equation}
     Therefore, in terms of the primitive variables the first subsystem reads $\partial_t \V + \hat{\mathbf{C}}_1 \partial_x \V = 0$, where
     \begin{equation}
    \label{eq:C_1system}
    \renewcommand{\arraystretch}{1.3}
         \hat{\mathbf{C}}_1 = \frac{\partial \V}{\partial \Q} \, \hat{\mathbf{A}}_1 \, \frac{\partial \Q}{\partial \V} = \begin{pmatrix}
         u & 0 & 0 & 0 & 0\\
         \dfrac{\rho_1}{\alpha_1}(u_1-u) & u_1 & \rho_1 & 0 & 0\\
         \dfrac{p_1-u_2^2 \rho_2}{\rho} & \dfrac{c_1^2}{\rho_1} & u_1 & \dfrac{u_2^2 \alpha_2}{\rho} & \dfrac{u_2 \alpha_2 \rho_2}{\rho}\\
         \dfrac{\rho_2 u}{\alpha_2} & 0 & 0 & 0 & 0 \\
         \dfrac{p_1-u_2^2 \rho_2}{\rho} & 0 & 0 & \dfrac{u_2^2 \alpha_2}{\rho} & \dfrac{u_2 \alpha_2 \rho_2}{\rho}+ u_2
         \end{pmatrix}
     \end{equation}
     The eigenvalues are then given by $\lambda_{1,1}, \lambda_{1,2} = u_1 \pm c_1$, $\lambda_{1,3} = u$, $\lambda_{1,4} = u_2\left(1+\frac{\alpha_2 \rho_2}{\rho} \right)$ and $\lambda_{1,5} = 0$.
     
     Let us investigate the eigenvectors of the system. For the acoustic waves $\lambda_{1,1}, \lambda_{1,2} = u_1 \pm c_1$ we get the following eigenvectors
    \begin{equation}
        \mathbf{r}_1 = \begin{pmatrix} 0 \\ \rho_1 \\ c_1 \\ 0 \\ 0 \end{pmatrix}, \quad 
        \mathbf{r}_2 = \begin{pmatrix} 0 \\ \rho_1 \\ -c_1 \\ 0 \\ 0 \end{pmatrix}
    \end{equation}
    
    For the interfacial wave $\lambda_{1,3} = u$, defining $D = (u_1-u)^2 - c_1^2$ and assuming $\lambda_{1,4} \neq u$ and $D \neq 0$, the corresponding eigenvector is
    \begin{equation}
    \renewcommand{\arraystretch}{1.3}
        \mathbf{r}_3 = 
        \begin{pmatrix} 
        1 \\ 
        \dfrac{\rho_1 K - \hat{C}_{21}(u_1-u)}{D} \\ 
        \dfrac{\hat{C}_{21} (c_1^2/\rho_1) - K (u_1-u)}{D} \\ 
        \dfrac{\rho_2}{\alpha_2} \\ 
        -\dfrac{p_1}{\rho(\lambda_{1,4} - u)} 
        \end{pmatrix},
    \end{equation}
    where $\hat{C}_{21} = \frac{\rho_1}{\alpha_1}(u_1-u)$ and $\hat{C}_{35} = \frac{u_2 \alpha_2 \rho_2}{\rho}$ are entries of $\hat{\mathbf{C}}_1$ in \eqref{eq:C_1system} and $K = \frac{p_1}{\rho} - \frac{p_1 \hat{C}_{35}}{\rho(\lambda_{1,4} - u)}$.
    
    For the wave associated with phase 2, $\lambda_{1,4} = u_2\left(1+\frac{\alpha_2 \rho_2}{\rho} \right)$, we obtain
    \begin{equation}
        \mathbf{r}_4 = \begin{pmatrix} 0 \\ \rho_1 u_2 \alpha_2 \rho_2 \\ (\lambda_{1,4} - u_1) u_2 \alpha_2 \rho_2 \\ 0 \\ \rho [ (u_1 - \lambda_{1,4})^2 - c_1^2 ] \end{pmatrix}
    \end{equation}
    
    Finally, for the stationary wave $\lambda_{1,5} = 0$, assuming $u_2 \neq 0$, the corresponding eigenvector is
    \begin{equation}
        \mathbf{r}_5 = \begin{pmatrix} 
        0 \\ -\rho_1 u_2^3 \alpha_2 \\ u_1 u_2^3 \alpha_2 \\ u_2 (c_1^2 - u_1^2)(\rho + \alpha_2 \rho_2) \\ -u_2^2 \alpha_2 (c_1^2 - u_1^2)
        \end{pmatrix}
    \end{equation}
    
    Except for $u = \lambda_{1,4}$ and for the coincidence of $u_1 \pm c_1$ with $u$, $0$ or $\lambda_{1,4}$, these five eigenvectors (with $\mathbf{r}_5$ replaced by the fourth unit vector for $u_2 = 0$) are linearly independent, so that the system is strongly hyperbolic. For $u_2 \neq 0$, the system is in general only weakly hyperbolic in the excluded cases, although a complete set of eigenvectors can still exist at particular states, \textit{e.g.} for $(u_1-u)^2 = c_1^2$ and $K = c_1^2/\alpha_1$. For $u_2 = 0$, the matrix \eqref{eq:C_1system} coincides with \eqref{eq:C_gas}, so that Proposition~\ref{prop:gas_system} applies.
   
\end{proof}

\begin{rem}
    \label{rem:lambda14_bound}
    Let us observe that, since $\alpha_2 \rho_2 \leq \rho$ by definition, one has
    \begin{equation}
        |\lambda_{1,4}| \leq 2|u_2|.
    \end{equation}
\end{rem}
\subsubsection{The second subsystem}
The second subsystem is given by
\begin{subequations}\label{eq:2-system}
\begin{eqnarray}
    \partial_t \rho + \partial_x (\alpha_2 \rho_2 u_2) &=& 0, \\
    \partial_t \alpha_1 &=& -\frac{1}{\tau}(p_2-p_1), \\
    \partial_t (\alpha_1 \rho_1) &=& 0, \\
    \partial_t (\rho u) + \partial_x (\alpha_2 p_2) &=& 0, \\
    \partial_t (u_1-u_2) - \partial_x (h_2) &=& -\zeta \frac{\alpha_1 \rho_1 \alpha_2 \rho_2}{\rho}(u_1 - u_2).
\end{eqnarray}
\end{subequations}

The second subsystem \eqref{eq:2-system} retains the triple eigenvalue $0$, while the acoustic eigenvalues become
\begin{equation}
    \lambda_{4,5} = \frac{1}{2 \rho}\left(\alpha_1 \rho_1 u_2 \pm \sqrt{\alpha_1^2 \rho_1^2 u_2^2 + 4 \rho^2 c_2^2}\right).
\end{equation}
Since the radicand is a sum of squares, they are real and distinct. As for the liquid-system, however, the eigenvalue $0$ has only two linearly independent eigenvectors for $p_2 \neq 0$, so that the second subsystem is weakly hyperbolic.

\section{The Semi-Implicit Schemes}
\label{Section:semi_implicit}
In this section, we describe the numerical schemes for the two-phase flow model \eqref{system:two_phase}, which are based on a finite volume discretization. We present two semi-implicit schemes based on the two splitting strategies proposed in Section \ref{Section:flux_splitting}. In this work, the source terms $\mathbf{S}(\Q)$ are not taken into account in the numerical schemes, \textit{i.e.} only the homogeneous part of system \eqref{system:two_phase} is solved, which corresponds to the limit $\tau \to \infty$ and $\zeta = 0$. 
A classical explicit method is subject to a restriction on the time step of the form
\begin{equation}
    \Delta t_{\text{expl}} \leq \nu_{\text{CFL}} \frac{\Delta x}{\max(|u|, |u_1 \pm c_1|, |u_2 \pm c_2|)},
\end{equation}
which tends to zero together with the Mach number of the liquid phase. This motivates an implicit-explicit approach, which allows a time step that is independent of the fast waves of the stiff phase.

The scheme based on phase-splitting consists of an explicit stage for the gas-system \eqref{eq:gas_system} and an implicit stage for the liquid-system \eqref{eq:liquid-system}. This leads to a CFL condition of the form
\begin{equation}
    \Delta t_{\text{phase}} \leq \nu_{\text{CFL}} \frac{\Delta x}{\max(|\lambda^u|, |\lambda_1^\pm|)} = \nu_{\text{CFL}} \frac{\Delta x}{\max(|u|, |u_1 \pm c_1|)},
\end{equation}
\textit{i.e.} the time step follows the material wave $\lambda^u$ and the gas phase, while the liquid phase is not resolved in time.

The second approach, based on the alternative splitting, consists of an explicit stage for the first subsystem \eqref{eq:1-system} and an implicit stage for the second subsystem \eqref{eq:2-system}. In this case, the CFL condition becomes
\begin{equation}
    \label{eq:cfl_alternative}
    \Delta t_{\text{altern}} \leq \nu_{\text{CFL}} \frac{\Delta x}{\max(|u|, |u_1 \pm c_1|, |\lambda_{1,4}|)},
\end{equation}
\textit{i.e.} the time step follows the material wave $\lambda^u$, the gas phase and the convective velocity $|\lambda_{1,4}|$ of the liquid bulk, and only the acoustic waves of the liquid phase are not resolved in time. Since $|\lambda_{1,4}| \leq 2|u_2|$, see Remark~\ref{rem:lambda14_bound}, replacing $|\lambda_{1,4}|$ by $2|u_2|$ gives a simpler but more restrictive sufficient condition.

One of the main difficulties in the numerical approximation of two-phase flow models is the pressure--velocity preservation condition (also referred to as the \textit{Abgrall condition}), see \cite{abgrall1996, saurel1999, dumbser2010}. It requires that a state in which both phases move with the same uniform velocity and have the same uniform pressure, for an arbitrary volume fraction, namely
\begin{equation}
    \label{eq:Abgrall_condition}
    u_1 = u_2 = u = \text{const}, \qquad p_1 = p_2 = p = \text{const},
\end{equation}
is preserved exactly by the numerical scheme.

\subsection{Phase-Splitting Scheme}
\label{Section:numerical_phase_splitting}
\subsubsection{Explicit Stage}
The first stage employs a first-order finite volume scheme to solve the gas-system explicitly in the form
\begin{equation}
    \label{eq:explicit_phase}
    \Q^{\ast}_i = \Q^n_i - \frac{\Delta t}{\Delta x} \left( \mathcal{F}_{i+\halb}^n - \mathcal{F}_{i-\halb}^n \right) - \frac{\Delta t}{\Delta x} \left( \mathcal{B}_{i+\halb}^n + \mathcal{B}_{i-\halb}^n \right),
\end{equation}
with $\Delta x = x_{i+\halb} - x_{i-\halb}$ assumed to be constant for all cells and $\Delta t = t^{n+1}-t^n$.
For the conservative part we use the Rusanov flux, defined as 
\begin{equation}
    \label{eq:rusanov_phase}
    \mathcal{F}^n_{i+\halb} = \frac{1}{2}(\f_1(\Q_{i+1}^n) + \f_1(\Q_i^n)) - \frac{1}{2} s^{\max}_{i+\halb} \, \left( \Q_{i+1}^n - \Q_i^n \right),
\end{equation}
where $\f_1$ is the conservative physical flux of the gas-system defined in \eqref{eq:fluxes_phase} and 
\begin{equation}
    \label{eq:smax_phase}
    s^{\max}_{i+\halb} = \max\left(\max_k |\lambda_{1,k}(\Q^n_{i+1})|, \max_k |\lambda_{1,k}(\Q^n_{i})|\right).
\end{equation}

To preserve the \textit{Abgrall condition} \eqref{eq:Abgrall_condition}, the non-conservative part is approximated by
\begin{equation}
    \label{eq:ncp_phase}
    \mathcal{B}^n_{i+\halb} = \Biggl(0, \frac{1}{2} \bar{u}_{i+\halb} (\alpha^n_{1,i+1} - \alpha^n_{1,i}), 0, 0, 0 \Biggr)^\intercal
\end{equation}
with the average velocity 
\begin{equation}
    \label{eq:ubar_phase}
    \bar{u}_{i+\halb} = \frac{1}{2} (u_i^n+u_{i+1}^n).
\end{equation}
In the framework of path-conservative schemes \cite{castro2006,pares2006}, this average corresponds to the integration of the non-conservative product along the segment path connecting $\Q_i^n$ and $\Q_{i+1}^n$ with the trapezoidal rule.
As a result of this explicit stage, we obtain
\begin{equation}
    \Q_i^\ast = \begin{pmatrix}
    \rho_i^\ast\\
    \alpha_{1,i}^\ast\\
    (\alpha_1 \rho_1)_i^\ast\\
    (\rho u)_i^\ast\\
    (u_1-u_2)_i^\ast
    \end{pmatrix},
\end{equation}
where it is worth noting that $\alpha_{1,i}^{n+1}=\alpha_{1,i}^\ast$ and $(\rho_1)_i^{n+1}=(\rho_1)_i^\ast$, since the volume fraction and the gas density are no longer updated in the implicit stage (liquid-system).

\subsubsection{Implicit Stage}
We proceed with temporal semi-discretization, keeping space still continuous. Let us consider the equation for the mixture density update
\begin{equation}
\label{eq:phase_splitting_implicit_stage_density}
    \rho^{n+1} = \rho^\ast - \Delta t \partial_x (\alpha_2^{n+1} \rho_2^{n+1} u_2^{n+1}).
\end{equation}
Let us also consider the equation for the momentum, where $(\rho u)^\ast$ is the solution of the explicit stage,
\begin{equation}
    \label{eq:momentum_phase}
    (\rho u)^{n+1} = (\rho u)^\ast - \Delta t \partial_x (\alpha_2^{n} \rho_2^{n+1} (u_2^{n+1})^2 + \alpha_2^{n+1} p_2^{n+1}).
\end{equation}

Since Equation \eqref{eq:phase_splitting_implicit_stage_density} holds, we can rewrite
\begin{equation}
\label{eq:phase_splitting_implicit_stage_mom2}
    (\alpha_2 \rho_2 u_2)^{n+1} = (\rho u)^\ast - (\alpha_1 \rho_1 u_1)^{n+1} - \Delta t \partial_x \left(\alpha_2^{n} \rho_2^{n+1} (u_2^{n+1})^2 + \alpha_2^{n+1} p_2^{n+1} \right)
\end{equation}

We need to derive the equation for $u_1^{n+1}$ in terms of $u_2^{n+1}$. Let us consider the equation for the relative velocity, \textit{i.e.} the last equation of \eqref{eq:liquid-system}, which is discretized as
\begin{equation}
    \label{eq:u1_phase}
    u_1^{n+1} = u_2^{n+1} + (u_1-u_2)^\ast + \Delta t \partial_x \left(\frac{1}{2} (u_2^{n+1})^2 + h_2^{n+1} \right)
\end{equation}
Substituting this into \eqref{eq:phase_splitting_implicit_stage_mom2}, we obtain
\begin{multline}
    \label{eq:phase_splitting_implicit_stage_mom2_final} 
    (\alpha_2 \rho_2 u_2)^{n+1} = (\rho u)^\ast - (\alpha_1 \rho_1)^{n+1} \left[u_2^{n+1} + (u_1-u_2)^\ast + \Delta t \partial_x \left(\frac{1}{2} (u_2^{n+1})^2 + h_2^{n+1} \right) \right] \\ - \Delta t \partial_x \left(\alpha_2^{n} \rho_2^{n+1} (u_2^{n+1})^2 + \alpha_2^{n+1} p_2^{n+1} \right)
\end{multline}
Finally, inserting \eqref{eq:phase_splitting_implicit_stage_mom2_final} into \eqref{eq:phase_splitting_implicit_stage_density}, we get
\begin{multline}
    \rho^{n+1} = \rho^\ast - \Delta t \partial_x \Biggl[ (\rho u)^\ast - (\alpha_1 \rho_1)^{n+1} \Bigl(u_2^{n+1} + (u_1-u_2)^\ast \\ + \Delta t \partial_x \left(\frac{1}{2} (u_2^{n+1})^2 + h_2^{n+1} \right) \Bigr) - \Delta t \partial_x \left(\alpha_2^{n} \rho_2^{n+1} (u_2^{n+1})^2 + \alpha_2^{n+1} p_2^{n+1} \right)\Biggr],
\end{multline}
which leads to the following discrete second-order wave equation 
\begin{multline}
\label{eq:wave_phase}
    (\alpha_2 \rho_2)^{n+1}  - \Delta t^2 \partial_x^2 (\alpha_2^{n+1} p_2^{n+1}) - \Delta t^2 \partial_x^2 \left( \alpha_2^{n} \rho_2^{n+1} (u_2^{n+1})^2 \right) \\ = \rho^\ast - (\alpha_1 \rho_1)^{n+1} - \Delta t \partial_x \left[ (\rho u)^\ast - (\alpha_1 \rho_1 u_1)^{n+1} \right].
\end{multline}

The fully discrete phase-splitting scheme uses the same discretization as the alternative splitting scheme of Section~\ref{sec:alternative_splitting}, namely the momentum predictor \eqref{eq:momentum_corrected}, in which the explicit gas pressure contribution is replaced by the gas pressure at the new volume fraction, the liquid momentum \eqref{eq:mu_star} and the liquid mass flux \eqref{eq:mass_flux_implicit}--\eqref{eq:interface_densities}. For a given liquid velocity $u_2^{(k)}$, the interface momentum \eqref{eq:mu_interface} is extended by the lagged convective terms of the liquid,
\begin{eqnarray}
\label{eq:mu_phase}
    \mu_{i+\halb} &=& \overline{\mu}^*_{i+\halb} - \Delta t D^+_{x,i+\halb}\left(\alpha_1^{n+1} p_1^* + \alpha_2^{n+1} p_2\right) - \Delta t\, \overline{D_x\left(\alpha_2^n \rho_2 (u_2^{(k)})^2\right)}_{i+\halb} \nonumber \\
    & & - \, \Delta t\, \overline{(\alpha_1 \rho_1)^*}_{i+\halb} D^+_{x,i+\halb}\left(\frac{1}{2}(u_2^{(k)})^2 + h_2\right),
\end{eqnarray}
where the old volume fraction $\alpha_2^n$ in the convective term is necessary for the preservation of the Abgrall condition, see Section~\ref{sec:abgrall_condition}. Here $\overline{(\cdot)}_{i+\halb}$ denotes an average of the values in the two adjacent cells and $D^+_{x, i+\halb}f = (f_{i+1}-f_i)/\Delta x$ the forward difference operator. With the mass flux \eqref{eq:mass_flux_implicit} evaluated with the interface momentum \eqref{eq:mu_phase}, the discrete system for the liquid pressure takes the form
\begin{equation}
\label{eq:root_pressure}
    \mathcal{R}_i := \alpha_{2,i}^{n+1} \rho_2(p_{2,i}) - \left(\rho_i^\ast - (\alpha_1 \rho_1)_i^\ast\right) + \frac{\Delta t}{\Delta x} \left( \mathcal{F}^{\alpha_2\rho_2}_{i+\halb} - \mathcal{F}^{\alpha_2\rho_2}_{i-\halb} \right) = 0.
\end{equation}
Since the convective term is differenced twice, the residual couples five neighbouring cells. Once $p_2$ is known, the momentum and the relative velocity are updated as
\begin{eqnarray}
    (\rho u)_i &=& \widetilde{\rho u}^*_i - \Delta t\, D_x\left(\alpha_1^{n+1} p_1^* + \alpha_2^{n+1} p_2\right)_i - \Delta t\, D_x\left(\alpha_2^n \rho_2 (u_2^{(k)})^2\right)_i, \label{eq:update_phase} \\
    (u_1-u_2)_i &=& (u_1-u_2)^*_i + \Delta t\, D_x\left(\frac{1}{2}(u_2^{(k)})^2 + h_2\right)_i, \label{eq:update_phase_b}
\end{eqnarray}
and the liquid velocity is recovered by \eqref{eq:velocities_alternative}.

\begin{rem}
\label{remark_newtonraphson}
The residuals \eqref{eq:root_pressure} and \eqref{eq:discrete_residual_implicit} are solved for the liquid pressure by the Newton method, since the liquid density follows from the inverted equation of state $\rho_2 = \rho_2(p_2)$. All quantities in the residuals, including $\rho_2$, $h_2$ and the interface densities, are evaluated at the current Newton iterate, so that also the off-diagonal entries of the Jacobian depend on the pressure. The Jacobian is tridiagonal for the alternative splitting and pentadiagonal for the phase splitting, because of the centred convective term, and it is computed by finite differences, which require only three or five evaluations of the residual, respectively, if the perturbations are applied simultaneously to cells that do not belong to a common stencil. The initial guess is $p_2^n$, each Newton increment is limited such that $p_2 + p_\infty$ remains positive, and the iterations stop when $\| (p_2^{(k+1)} - p_2^{(k)}) / (p_2^{(k+1)} + p_\infty) \|_\infty < 10^{-12}$.
\end{rem}

The lagged liquid velocity $u_2^{(k)}$ is updated by Picard iterations, which are combined with the Newton method for the pressure: each iteration consists of one Newton step for the residual \eqref{eq:root_pressure}, followed by the update of $u_2$. The complete time step of the phase-splitting scheme can be summarized as follows.

\begin{algo}
\label{alg:picard_newton}
\mbox{}
\begin{enumerate}
    \item Compute $\Q_i^\ast$ by the explicit stage \eqref{eq:explicit_phase} with the Rusanov flux \eqref{eq:rusanov_phase}--\eqref{eq:smax_phase} and the non-conservative term \eqref{eq:ncp_phase}--\eqref{eq:ubar_phase}, set $\alpha_{1,i}^{n+1} = \alpha_{1,i}^\ast$ and $\rho_{1,i}^{n+1} = \rho_{1,i}^\ast$, and compute $p_{1,i}^\ast = p_1(\rho_{1,i}^\ast)$.
    \item Compute the momentum predictor $\widetilde{\rho u}^*_i$ and the liquid momentum $\mu_i^*$ by \eqref{eq:momentum_corrected}--\eqref{eq:mu_star}, and set $p_2^{(0)} = p_2^n$ and $u_2^{(0)} = u_2^n$.
    \item Compute $p_2^{(k+1)}$ by one Newton step for the residual \eqref{eq:root_pressure} with the interface momentum \eqref{eq:mu_phase}.
    \item Update $(\rho u)_i$ and $(u_1-u_2)_i$ by \eqref{eq:update_phase}--\eqref{eq:update_phase_b} with $p_2^{(k+1)}$, and recover $u_2^{(k+1)}$ by \eqref{eq:velocities_alternative}.
    \item If the relative increment of $p_2 + p_\infty$ and the increment $\|u_2^{(k+1)} - u_2^{(k)}\|_\infty / (\|u_2^{(k+1)}\|_\infty + \lambda_{\max})$ of the liquid velocity, where $\lambda_{\max} = \max_i \max(|u_i|, |u_{1,i}| + c_{1,i})$ is the largest wave speed of the explicit stage, are below $10^{-12}$, go to step 6; otherwise, set $k \leftarrow k+1$ and go back to step 3.
    \item Set $p_2^{n+1} = p_2^{(k+1)}$ and $\rho_2^{n+1} = \rho_2(p_2^{n+1})$, and compute $\rho^{n+1} = (\alpha_1 \rho_1)^\ast + \alpha_2^{n+1} \rho_2^{n+1}$, $(\rho u)^{n+1}$ and $(u_1-u_2)^{n+1}$ by \eqref{eq:update_phase}--\eqref{eq:update_phase_b} with $u_2^{(k+1)}$, and $u_2^{n+1}$ and $u_1^{n+1}$ by \eqref{eq:velocities_alternative}--\eqref{eq:velocities_alternative_b}.
\end{enumerate}
\end{algo}

\subsection{Alternative Splitting Scheme}\label{sec:alternative_splitting}
The main advantage of the alternative splitting scheme presented in this section is that a nonlinear discrete wave equation for the pressure can be solved directly, without nested Picard iterations in the velocity, while maintaining a moderate CFL condition. This is achieved by treating the entire gasdynamics and the material advection terms of the liquid phase fully explicitly, so that only the fast acoustic waves of the liquid phase are treated implicitly. As for the phase splitting, the algorithm is split into two stages: an explicit prediction step featuring a consistency correction, followed by an implicit correction step.

\subsubsection{Explicit Stage}
The first stage of the algorithm solves the first subsystem, which governs the slow material waves and the gasdynamics, explicitly by a first-order finite volume scheme
\begin{equation}
\label{eq:finite_volume_explicit}
    \mathbf{Q}^{\ast}_i = \mathbf{Q}^n_i - \frac{\Delta t}{\Delta x} \left( \mathcal{F}_{i+\halb}^n - \mathcal{F}_{i-\halb}^n \right) - \frac{\Delta t}{\Delta x} \left( \mathcal{B}_{i+\halb}^n + \mathcal{B}_{i-\halb}^n \right),
\end{equation}
where $\mathbf{Q} = [\rho, \alpha_1, \alpha_1\rho_1, \rho u, u_1-u_2]^T$. The numerical fluxes at the cell interfaces are again given by the local Rusanov flux
\begin{equation}
    \label{eq:rusanov_alternative}
    \mathcal{F}^n_{i+\halb} = \frac{1}{2}(\hat{\mathbf{f}}_1(\mathbf{Q}_{i+1}^n) + \hat{\mathbf{f}}_1(\mathbf{Q}_i^n)) - \frac{1}{2} s^{\max}_{i+\halb} \, \left( \mathbf{Q}_{i+1}^n - \mathbf{Q}_i^n \right),
\end{equation}
where the maximum local wave speed is selected according to the CFL condition
\begin{equation}
    \label{eq:smax_alternative}
    s^{\max}_{i+\halb} = \max_{k=\{i,i+1\}} \left( |u_k|, |u_{1,k} \pm c_{1,k}^n|, |\lambda_{1,4}(\mathbf{Q}^n_k)| \right),
    \qquad \lambda_{1,4} = u_2\left(1+\frac{\alpha_2\rho_2}{\rho}\right).
\end{equation}

To satisfy the \textit{Abgrall condition} \eqref{eq:Abgrall_condition} and to avoid spurious pressure oscillations across moving material fronts, the non-conservative advection term of the volume fraction is again discretized symmetrically with the interface velocities $\bar{u}_{i+\halb} = \frac{1}{2}(u_i^n + u_{i+1}^n)$, namely
\begin{equation}
\label{eq:approximation_NC_term}
     \mathcal{B}^n_{i+\halb} = \Biggl(0, \frac{1}{2} \bar{u}_{i+\halb} (\alpha^n_{1,i+1} - \alpha^n_{1,i}), 0, 0, 0 \Biggr)^\intercal.
\end{equation}
As a result of this explicit stage, we obtain
\begin{equation}
    \Q_i^\ast = \begin{pmatrix}
    \rho_i^\ast\\
    \alpha_{1,i}^\ast\\
    (\alpha_1 \rho_1)_i^\ast\\
    (\rho u)_i^\ast\\
    (u_1-u_2)_i^\ast
    \end{pmatrix},
\end{equation}
where again it is worth noting that $\alpha_{1,i}^{n+1}=\alpha_{1,i}^\ast$ and $(\rho_1)_i^{n+1}=(\rho_1)_i^\ast$.

Following this explicit update, the volume fraction is advanced to the next time level by setting $\alpha_1^{n+1} = \alpha_1^*$ and $\alpha_2^{n+1} = 1 - \alpha_1^{n+1}$. The updated gas density and the explicit gas pressure are then computed as 
\begin{equation}
    \label{eq:rho1_alternative}
\rho_1^* = (\alpha_1\rho_1)^*/\alpha_1^{n+1}, \qquad p_1^* = p_1(\rho_1^*).
\end{equation}

Finally, the liquid mass is computed as $(\alpha_2\rho_2)^* = \rho^* - (\alpha_1\rho_1)^*$.

\subsubsection{Implicit Correction and Velocity Recovery}
The second stage updates the liquid phase mass by accounting for the fast acoustic waves. Let us consider the equation for the mixture density update in the implicit stage
\begin{equation}
    \rho^{n+1} = \rho^\ast - \Delta t \partial_x (\alpha_2^{n+1} \rho_2^{n+1} u_2^{n+1}).
\end{equation}
Let us consider the equation for the momentum. Since the convective part of the liquid phase is now treated explicitly, the implicit update only contains the liquid pressure gradient
\begin{equation}
    (\rho u)^{n+1} = (\rho u)^\ast - \Delta t \partial_x (\alpha_2^{n+1} p_2^{n+1}).
\end{equation}
Using $(\rho u) = \alpha_1 \rho_1 u_1 + \alpha_2 \rho_2 u_2$, we can rewrite the liquid momentum as
\begin{equation}
\label{alternative_momentum_discretized}
    (\alpha_2 \rho_2 u_2)^{n+1} = (\rho u)^\ast - (\alpha_1 \rho_1 u_1)^{n+1} - \Delta t \partial_x \left(\alpha_2^{n+1} p_2^{n+1} \right).
\end{equation}
We also consider the equation for the relative velocity, \textit{i.e.} the last equation of \eqref{eq:2-system}, discretized as
\begin{equation}
\label{alternative_relative_discretized}
    u_1^{n+1} = u_2^{n+1} + (u_1-u_2)^\ast + \Delta t \partial_x (h_2^{n+1}).
\end{equation}
From \eqref{alternative_momentum_discretized} and \eqref{alternative_relative_discretized}, we finally obtain
\begin{equation}
    \rho^{n+1} u_2^{n+1} = (\rho u)^\ast - (\alpha_1 \rho_1)^\ast (u_1 - u_2)^\ast - \Delta t (\alpha_1 \rho_1)^\ast \partial_x (h_2^{n+1}) - \Delta t \partial_x \left(\alpha_2^{n+1} p_2^{n+1} \right).
\end{equation}
Substituting this into the density equation, we get
\begin{multline}
    \rho^{n+1} = \rho^\ast - \Delta t \partial_x \Biggl[\frac{\alpha_2^{n+1} \rho_2^{n+1}}{\rho^{n+1}} \Biggl( (\rho u)^\ast - (\alpha_1 \rho_1)^\ast (u_1 - u_2)^\ast \\ - \Delta t (\alpha_1 \rho_1)^\ast \partial_x (h_2^{n+1}) - \Delta t \partial_x \left(\alpha_2^{n+1} p_2^{n+1} \right) \Biggr) \Biggr]
\end{multline}
Since $\rho = \alpha_1 \rho_1 + \alpha_2 \rho_2$ and the gas mass is frozen, $(\alpha_1 \rho_1)^{n+1} = (\alpha_1 \rho_1)^\ast$, we obtain the following discrete second-order wave equation
\begin{multline}
\label{eq:wave_alternative}
    (\alpha_2 \rho_2)^{n+1}  = \rho^\ast - (\alpha_1 \rho_1)^\ast - \Delta t \partial_x \Biggl[\frac{\alpha_2^{n+1} \rho_2^{n+1}}{\rho^{n+1}} \Biggl( (\rho u)^\ast - \\   (\alpha_1 \rho_1)^\ast (u_1 - u_2)^\ast - \Delta t (\alpha_1 \rho_1)^\ast \partial_x (h_2^{n+1}) - \Delta t \partial_x \left(\alpha_2^{n+1} p_2^{n+1} \right) \Biggr) \Biggr],
\end{multline}
where the implicit mixture density is $\rho^{n+1} = (\alpha_1 \rho_1)^\ast + (\alpha_2 \rho_2)^{n+1}$, and the liquid state variables $\rho_2^{n+1}$ and $h_2^{n+1}$ are expressed through the liquid EOS \eqref{eq:eos_p2} as functions of $p_2^{n+1}$.

To discretize \eqref{eq:wave_alternative}, we define for each cell $i$ the nonlinear residual
\begin{equation}
    \label{eq:discrete_residual_implicit}
    \mathcal{R}_i(p_2^{n+1}) := \alpha_{2,i}^{n+1} \rho_2(p_{2,i}^{n+1}) - \left(\rho_i^\ast - (\alpha_1 \rho_1)_i^\ast\right) + \frac{\Delta t}{\Delta x} \left( \mathcal{F}^{\alpha_2\rho_2}_{i+\halb} - \mathcal{F}^{\alpha_2\rho_2}_{i-\halb} \right) = 0,
\end{equation}
where the discrete numerical mass flux $\mathcal{F}^{\alpha_2\rho_2}_{i+\halb}$ is constructed by splitting it into a structure-preserving centered formulation and a dissipation term approximating the square bracketed expression.

The main difficulty of this approach is the preservation of the \textit{Abgrall condition} \eqref{eq:Abgrall_condition}. For a state of constant pressure $p_1 = p_2 = \bar{p}$, also $\alpha_1 p_1 + \alpha_2 p_2 = \bar{p}$ is constant. With a splitting in which $\alpha_1 p_1$ is treated explicitly and $\alpha_2 p_2$ implicitly, such states are in general not preserved exactly. Indeed, a direct evaluation of the mass flux with the updated volume fraction $\alpha_2^{n+1}$ generates spurious pressure oscillations across moving material fronts, due to a discrete mismatch between the explicit update of the gas volume fraction and the implicit update of the liquid phase.

To preserve these states, the implicit stage must use the old volume fraction $\alpha_2^n$ together with the new implicit density $\rho_2^{n+1}$. Let us define a modified, purely convective total momentum predictor $\widetilde{\rho u}^*_i$ by stripping out the explicit gas pressure contribution and replacing it with a non-dissipative centered gradient
\begin{equation}
\label{eq:momentum_corrected}
    \widetilde{\rho u}^*_i = (\rho u)^*_i + \Delta t D_x (\alpha_1^n p_1^n)_i,
\end{equation}
where $D_x f_i = (f_{i+1} - f_{i-1})/(2\Delta x)$ is the standard second-order centered gradient. Then, we can isolate the cell-centered liquid momentum driving the flux in \eqref{eq:discrete_residual_implicit}
\begin{equation}
    \label{eq:mu_star}
    \mu_i^* = \widetilde{\rho u}^*_i - (\alpha_1 \rho_1)_i^* (u_1-u_2)^*_i.
\end{equation} 
The implicitly corrected interface momentum, which corresponds to the expression in round brackets in \eqref{eq:wave_alternative}, is then
\begin{equation}
    \label{eq:mu_interface}
    \mu_{i+\halb}^{n+1} = \overline{\mu}^*_{i+\halb} - \Delta t D^+_{x, i+\halb}\left(\alpha_1^{n+1} p_1^* + \alpha_2^{n+1} p_2^{n+1}\right) - \Delta t \overline{(\alpha_1 \rho_1)^*}_{i+\halb} D^+_{x, i+\halb}\left(h_2^{n+1}\right),
\end{equation}
where $\overline{(\cdot)}_{i+\halb}$ denotes an average of the values in the two adjacent cells and $D^+_{x, i+\halb}f = (f_{i+1}-f_i)/\Delta x$ the forward difference operator.

To guarantee the preservation of the Abgrall condition, the implicit liquid mass flux uses the old volume fraction $\alpha_2^n$ together with the new implicit density $\rho_2^{n+1}$, consistently with the discrete balance of the explicit advection step, namely
\begin{equation}
    \label{eq:mass_flux_implicit}
    \mathcal{F}^{\alpha_2\rho_2}_{i+\halb} = \left( \frac{\left(\alpha_2^n\rho_2^{n+1}\right)_{i+\halb}}{\overline{\rho}^{n+1}_{i+\halb}} \right) \mu_{i+\halb}^{n+1},
\end{equation}
where the implicit mixture density and the liquid mass field at the interface are defined by
\begin{equation}
    \label{eq:interface_densities}
    \overline{\rho}_{i+\halb}^{n+1} = \frac{\rho_i^{n+1} + \rho_{i+1}^{n+1}}{2}, \qquad \left(\alpha_2^n\rho_2^{n+1}\right)_{i+\halb} = \frac{\alpha_{2,i}^n \rho_2(p_{2,i}^{n+1}) + \alpha_{2,i+1}^n \rho_2(p_{2,i+1}^{n+1})}{2}.
\end{equation}
This choice of the interface values ensures the pressure--velocity preservation, as shown in Section \ref{sec:abgrall_condition}.
The residual \eqref{eq:discrete_residual_implicit} is solved for $p_2^{n+1}$ by the Newton method described in Remark~\ref{remark_newtonraphson}. After convergence, the liquid densities $\rho_2^{n+1}$ and enthalpies $h_2^{n+1}$ are computed, and the momentum and the relative velocity are updated as
\begin{eqnarray}
    \label{eq:update_alternative} (\rho u)_i^{n+1} &=& \widetilde{\rho u}^*_i - \Delta t \left[ D_x \left( \alpha_1^{n+1} p_1^* + \alpha_2^{n+1} p_2^{n+1} \right) \right]_i, \\
    \label{eq:update_alternative_b} (u_1-u_2)_i^{n+1} &=& (u_1-u_2)^*_i + \Delta t \left[ D_x h_2^{n+1} \right]_i.
\end{eqnarray}
From the total momentum $(\rho u)^{n+1}$, the phase velocities are recovered as
\begin{eqnarray}
    \label{eq:velocities_alternative} u_{2,i}^{n+1} &=& \frac{(\rho u)_i^{n+1} - (\alpha_1\rho_1)_i^* (u_1-u_2)_i^{n+1}}{\rho_i^{n+1}}, \\
    \label{eq:velocities_alternative_b} u_{1,i}^{n+1} &=& u_{2,i}^{n+1} + (u_1-u_2)_i^{n+1}.
\end{eqnarray}

The complete time step of the alternative splitting scheme can be summarized as follows: 

\begin{algo}
\label{alg:alternative}
\mbox{}
\begin{enumerate}
    \item Compute $\Q_i^\ast$ by the explicit stage \eqref{eq:finite_volume_explicit} with the Rusanov flux \eqref{eq:rusanov_alternative}--\eqref{eq:smax_alternative} and the non-conservative term \eqref{eq:approximation_NC_term}, set $\alpha_{1,i}^{n+1} = \alpha_{1,i}^\ast$ and $\alpha_{2,i}^{n+1} = 1 - \alpha_{1,i}^{n+1}$, and compute $\rho_{1,i}^\ast$ and $p_{1,i}^\ast$ by \eqref{eq:rho1_alternative}.
    \item Compute the momentum predictor $\widetilde{\rho u}^*_i$ and the liquid momentum $\mu_i^\ast$ by \eqref{eq:momentum_corrected} and \eqref{eq:mu_star}.
    \item Compute $p_2^{n+1}$ by the Newton method applied to the residual \eqref{eq:discrete_residual_implicit}, with the interface momentum \eqref{eq:mu_interface} and the mass flux \eqref{eq:mass_flux_implicit}--\eqref{eq:interface_densities}.
    \item Compute $\rho_2^{n+1}$ and $h_2^{n+1}$ from the equation of state \eqref{eq:eos_p2} and $\rho^{n+1} = (\alpha_1 \rho_1)^\ast + \alpha_2^{n+1} \rho_2^{n+1}$, and update $(\rho u)^{n+1}$ and $(u_1-u_2)^{n+1}$ by \eqref{eq:update_alternative}--\eqref{eq:update_alternative_b}.
    \item Recover $u_2^{n+1}$ and $u_1^{n+1}$ by \eqref{eq:velocities_alternative}--\eqref{eq:velocities_alternative_b}.
\end{enumerate}
\end{algo}

\subsection{Pressure--Velocity Preservation Condition of the Schemes}
\label{sec:abgrall_condition}
To prove the pressure--velocity preservation property (\textit{Abgrall condition}) for material interfaces, we analyze the preservation of an initial state at time level $t^n$ characterized by a uniform velocity and a uniform pressure across both phases, namely
\begin{equation}
    u_{1,i}^n = u_{2,i}^n = \bar{u}, \quad p_{1,i}^n = p_{2,i}^n = \bar{p}, \quad \forall i,
\end{equation}
with boundary conditions that are compatible with this state.

Under isentropic equations of state, a uniform pressure implies constant thermodynamic phase densities, $\rho_{1,i}^n = \bar{\rho}_1$ and $\rho_{2,i}^n = \bar{\rho}_2$. The volume fractions $\alpha_{1,i}^n$ and $\alpha_{2,i}^n = 1 - \alpha_{1,i}^n$ vary arbitrarily in space. The total mixture density and momentum are given by $\rho_i^n = \alpha_{1,i}^n\bar{\rho}_1 + \alpha_{2,i}^n\bar{\rho}_2$ and $(\rho u)_i^n = \rho_i^n \bar{u}$, respectively.

In the explicit stage \eqref{eq:finite_volume_explicit}, the central part of the Rusanov flux \eqref{eq:rusanov_alternative} yields the centered difference $D_x$ of the physical flux, with $D_x$ as defined after \eqref{eq:momentum_corrected}, and for a uniform velocity the non-conservative term \eqref{eq:approximation_NC_term} yields $\bar{u} D_x \alpha_{1,i}^n$. Since the same wave speed $s^{\max}_{i+\halb}$ is used for all components of $\Q$, the numerical viscosity adds the term $\Delta t\, \mathcal{V}_i(q^n)$ to the update of each component $q$ of $\Q$, where
\begin{equation}
    \label{eq:viscosity_operator}
    \mathcal{V}_i(q) = \frac{1}{2 \Delta x} \left[ s^{\max}_{i+\halb} \left( q_{i+1} - q_i \right) - s^{\max}_{i-\halb} \left( q_i - q_{i-1} \right) \right]
\end{equation}
is linear in $q$ and vanishes for constant $q$. The volume fractions are therefore updated as
\begin{eqnarray}
    \alpha_{1,i}^* &=& \alpha_{1,i}^n - \Delta t \bar{u} D_x \alpha_{1,i}^n + \Delta t\, \mathcal{V}_i(\alpha_1^n), \label{eq:alpha_up} \\
    \alpha_{2,i}^{n+1} &=& 1 - \alpha_{1,i}^* = \alpha_{2,i}^n - \Delta t \bar{u} D_x \alpha_{2,i}^n + \Delta t\, \mathcal{V}_i(\alpha_2^n). \label{eq:alpha2_up}
\end{eqnarray}
The gas partial density is updated with the physical flux $\alpha_1 \rho_1 u_1 = \bar{u} \bar{\rho}_1 \alpha_1$, which yields
\begin{equation}
    (\alpha_1 \rho_1)_i^* = \bar{\rho}_1 \left( \alpha_{1,i}^n - \Delta t \bar{u} D_x \alpha_{1,i}^n + \Delta t\, \mathcal{V}_i(\alpha_1^n) \right) = \bar{\rho}_1 \alpha_{1,i}^*.
\end{equation}
As a consequence, the intermediate gas density and pressure are preserved exactly,
\begin{equation}
    \rho_{1,i}^* = \frac{(\alpha_1 \rho_1)_i^*}{\alpha_{1,i}^*} = \bar{\rho}_1 \implies p_{1,i}^* = p_1(\rho_{1,i}^*) = \bar{p}.
\end{equation}

In this splitting strategy, the explicit mixture density $\rho_i^*$ is updated through the gas phase mass flux. Therefore,
\begin{equation}
    \rho_i^* = \rho_i^n - \Delta t \bar{u} \bar{\rho}_1 D_x \alpha_{1,i}^n + \Delta t\, \mathcal{V}_i(\rho^n),
\end{equation}
and, by the linearity of $\mathcal{V}_i$, the intermediate liquid partial mass becomes
\begin{eqnarray}\label{eq:liq_mass_explicit}
    (\alpha_2 \rho_2)_i^* &=& \rho_i^* - (\alpha_1 \rho_1)_i^* = \rho_i^n - \bar{\rho}_1 \alpha_{1,i}^n + \Delta t\, \mathcal{V}_i(\rho^n - \bar{\rho}_1 \alpha_1^n) \nonumber \\
    &=& \bar{\rho}_2 \left( \alpha_{2,i}^n + \Delta t\, \mathcal{V}_i(\alpha_2^n) \right).
\end{eqnarray}
Hence, the liquid partial mass is not transported in the explicit stage, but it is modified by the numerical viscosity. The explicit update of the mixture momentum contains the convective flux $\bar{u}^2 \rho$, the gas pressure flux $\alpha_1 \bar{p}$ and the numerical viscosity, so that the corrected momentum \eqref{eq:momentum_corrected} becomes
\begin{eqnarray}\label{eq:rhou_c}
    \widetilde{\rho u}^*_i &=& (\rho u)_i^* + \Delta t D_x(\alpha_{1,i}^n \bar{p}) = \bar{u} \left( \rho_i^n - \Delta t \bar{u} D_x \rho_i^n + \Delta t\, \mathcal{V}_i(\rho^n) \right) \nonumber \\
    &=& \bar{u} \left( \bar{\rho}_1 \alpha_{1,i}^* + \bar{\rho}_2 \alpha_{2,i}^{n+1} \right),
\end{eqnarray}
where we have used \eqref{eq:alpha_up}--\eqref{eq:alpha2_up}. Since the flux of the relative velocity is constant and $(u_1-u_2)_i^n = 0$, also $(u_1-u_2)_i^* = 0$.

During the implicit stage, we must verify that the uniform pressure $p_{2,i}^{n+1} = \bar{p}$ is a root of the residual \eqref{eq:discrete_residual_implicit}, \textit{i.e.} $\mathcal{R}_i(\bar{p}) = 0$. For $p_{2,i}^{n+1} = \bar{p}$, the new mixture density is $\rho_i^{n+1} = (\alpha_1 \rho_1)_i^* + \alpha_{2,i}^{n+1} \bar{\rho}_2 = \bar{\rho}_1 \alpha_{1,i}^* + \bar{\rho}_2 \alpha_{2,i}^{n+1}$, so that \eqref{eq:mu_star} and \eqref{eq:rhou_c} give $\mu_i^* = \widetilde{\rho u}^*_i = \bar{u} \rho_i^{n+1}$. Moreover, the mixture pressure $\alpha_{1,i}^{n+1} p_{1,i}^* + \alpha_{2,i}^{n+1} \bar{p} = \bar{p}$ is spatially uniform, and the liquid enthalpy $h_{2,i} = \bar{h}_2$ is constant. Consequently, the pressure and enthalpy differences in \eqref{eq:mu_interface} vanish, and the interface momentum simplifies to
\begin{equation}
    \mu_{i+\halb}^{n+1} = \overline{\mu}^*_{i+\halb} = \bar{u}\, \overline{\rho}^{n+1}_{i+\halb},
\end{equation}
with the interface density \eqref{eq:interface_densities}. The liquid mass flux \eqref{eq:mass_flux_implicit} and its difference are then
\begin{eqnarray}
    \mathcal{F}^{\alpha_2\rho_2}_{i+\halb} &=& \frac{\left(\alpha_2^n \bar{\rho}_2\right)_{i+\halb}}{\overline{\rho}^{n+1}_{i+\halb}}\, \bar{u}\, \overline{\rho}^{n+1}_{i+\halb} = \bar{u} \bar{\rho}_2 \, \frac{\alpha_{2,i}^n + \alpha_{2,i+1}^n}{2}, \\
    \frac{\mathcal{F}^{\alpha_2\rho_2}_{i+\halb} - \mathcal{F}^{\alpha_2\rho_2}_{i-\halb}}{\Delta x} &=& \bar{u} \bar{\rho}_2 D_x \alpha_{2,i}^n. \label{eq:flux_difference_equilibrium}
\end{eqnarray}
Substituting \eqref{eq:flux_difference_equilibrium} together with \eqref{eq:alpha2_up} and \eqref{eq:liq_mass_explicit} into the residual \eqref{eq:discrete_residual_implicit} yields
\begin{eqnarray}
    \mathcal{R}_i(\bar{p}) &=& \bar{\rho}_2 \left( \alpha_{2,i}^n - \Delta t \bar{u} D_x \alpha_{2,i}^n + \Delta t\, \mathcal{V}_i(\alpha_2^n) \right) - \bar{\rho}_2 \left( \alpha_{2,i}^n + \Delta t\, \mathcal{V}_i(\alpha_2^n) \right) \nonumber \\
    & & + \, \Delta t \bar{u} \bar{\rho}_2 D_x \alpha_{2,i}^n = 0.
\end{eqnarray}
The viscosity terms of $\alpha_{2,i}^{n+1}$ and of $(\alpha_2 \rho_2)_i^*$ cancel each other, and the transport terms cancel because the liquid mass flux \eqref{eq:mass_flux_implicit} uses the old volume fraction $\alpha_2^n$. Therefore, the constant pressure $p_{2,i}^{n+1} = \bar{p}$ is a root of the residual, and if the Newton method is initialized with $p_{2,i}^n = \bar{p}$, the residual vanishes already for the initial guess. Finally, in the velocity recovery stage, the vanishing pressure and enthalpy gradients in \eqref{eq:update_alternative}--\eqref{eq:update_alternative_b} imply that the relative velocity remains zero, $(u_1-u_2)_i^{n+1} = 0$, and that the mixture momentum is updated as $(\rho u)_i^{n+1} = \widetilde{\rho u}^*_i = \bar{u} \rho_i^{n+1}$. Reconstructing the phase velocities by \eqref{eq:velocities_alternative}--\eqref{eq:velocities_alternative_b} yields
\begin{equation}
    u_{2,i}^{n+1} = \frac{(\rho u)_i^{n+1} - (\alpha_1 \rho_1)_i^* (u_1-u_2)_i^{n+1}}{\rho_i^{n+1}} = \bar{u}, \quad u_{1,i}^{n+1} = u_{2,i}^{n+1} + (u_1-u_2)_i^{n+1} = \bar{u}.
\end{equation}
This concludes the proof that the scheme preserves uniform pressure and velocity across material interfaces and thus satisfies the Abgrall condition.

The same result holds for the phase-splitting scheme of Section~\ref{Section:numerical_phase_splitting}. There, the explicit stage does not contain the convective flux of the liquid, so that the momentum predictor is
\begin{equation}
    \widetilde{\rho u}^*_i = \bar{u} \left( \rho_i^n - \Delta t \bar{u} \bar{\rho}_1 D_x \alpha_{1,i}^n + \Delta t\, \mathcal{V}_i(\rho^n) \right),
\end{equation}
but for $u_2^{(k)} = \bar{u}$ the lagged convective term in \eqref{eq:mu_phase} and \eqref{eq:update_phase} contributes $-\Delta t \bar{u}^2 \bar{\rho}_2 D_x \alpha_{2,i}^n$. Since $\bar{\rho}_1 D_x \alpha_1^n + \bar{\rho}_2 D_x \alpha_2^n = D_x \rho^n$, the same interface momentum and the same residual as above are obtained, while the lagged kinetic energy is constant. Hence, the uniform state is a fixed point of the combined Newton and Picard iterations, which are initialized with $p_2^n = \bar{p}$ and $u_2^n = \bar{u}$.

\subsection{Consistency with the Low Mach Number Limit of the Liquid}
\label{sec:low_mach_consistency}
Let us assume that $M_1 = M_{\text{gas}} = \mathcal{O}(1)$ and $M_2 = M_{\text{liquid}} \ll 1$, namely the gas phase is compressible and the liquid phase is weakly compressible. Furthermore, let us assume well-prepared initial data in the sense of Lemma~\ref{lem:well_prepared}, \textit{i.e.} $\Q^n \in \Omega_1^{wp}$.

The goal is to show that, in the limit $M_2 \to 0$, \textit{i.e.} when the liquid becomes incompressible and its sound speed tends to infinity, both schemes are consistent with the volume balance of the limit system. For the well-prepared data, the pressure must then enforce the divergence-free constraint 
\begin{equation}
    \label{eq:AP_constraint}
    \partial_x (\alpha_2 u_2)^{n+1} = 0.
\end{equation}

The following argument is formal: it assumes smooth discrete solutions, bounded admissible pressures and explicit wave speeds, and volume fractions that stay away from zero and one. In both schemes, the explicit stage transports $\rho$ and $\alpha_1 \rho_1$ with the same convective flux, and the implicit stage updates $\rho$ with the liquid mass flux, while $\alpha_1 \rho_1$ is frozen. Hence both schemes satisfy the discrete liquid mass balance
\begin{equation}
\label{eq:AP_finalstep_implicit_alternative}
    \frac{(\alpha_2 \rho_2)^{n+1} - (\alpha_2 \rho_2)^\ast}{\Delta t} + \partial_x (\alpha_2 \rho_2 u_2)^{n+1} = 0,
\end{equation}
where $(\alpha_2 \rho_2)^\ast = \rho^\ast - (\alpha_1 \rho_1)^\ast$ differs from $(\alpha_2 \rho_2)^n$ only by the numerical dissipation of the explicit stage. Here, $\partial_x (\alpha_2 \rho_2 u_2)^{n+1}$ denotes the discrete divergence of the implicit liquid mass flux of the respective scheme. In both schemes, this is the flux \eqref{eq:mass_flux_implicit}, which uses $\alpha_2^n$ instead of $\alpha_2^{n+1}$ together with interface averages. Since $\alpha_2^{n+1} - \alpha_2^n = \mathcal{O}(\Delta t)$, this adds an error of order $\mathcal{O}(\Delta t)$ to the error of order $\mathcal{O}(\Delta x)$ caused by the averages. For bounded pressures, the liquid density becomes constant, $\rho_2 = \bar{\rho}_2$, in the limit $M_2 \to 0$, and \eqref{eq:AP_finalstep_implicit_alternative} becomes
\begin{equation}
\label{eq:AP_limit_volume}
    \frac{\alpha_2^{n+1} - \alpha_2^{n}}{\Delta t} + \partial_x (\alpha_2 u_2)^{n+1} = \mathcal{O}(\Delta t + \Delta x),
\end{equation}
\textit{i.e.} a consistent discretization of the volume balance $\partial_t \alpha_2 + \partial_x (\alpha_2 u_2) = 0$ of the limit system. If, in addition, the volume fraction remains stationary at leading order, as for the well-prepared data of Lemma~\ref{lem:well_prepared}, this reduces to the incompressibility constraint \eqref{eq:AP_constraint}. The liquid pressure $p_2^{n+1}$ is the variable that enforces this constraint through the implicit stage. We restrict ourselves here to this consistency with the limit constraint; the uniform stability of the schemes with respect to $M_2$ is not addressed.

\section{Numerical Results}
\label{Section:results}
In this section, we discuss the numerical results obtained with the two semi-implicit splitting schemes presented in Section \ref{Section:semi_implicit}. The initial data are listed in Table \ref{table:initial_condition} and the reference quantities in Table \ref{table:ref_quantities}. In all tests, the initial discontinuity is located at $x = L/2$.

\begin{table}[ht]
\centering
\footnotesize
\setlength{\tabcolsep}{4.5pt}
\begin{tabular}{lcccccccccc}
\hline
Test & $\rho_{1L}$ & $\rho_{1R}$ & $u_{1L}$ & $u_{1R}$ & $\rho_{2L}$ & $\rho_{2R}$ & $u_{2L}$ & $u_{2R}$ & $\alpha_{1L}$ & $\alpha_{1R}$ \\
\hline
1 & 1.2449 & 0.60312 & $-1.2638$ & 0.43059 & 1.2969 & 0.73436 & $-0.3894$ & $-0.40507$ & 0.7 & 0.3 \\
2a & 1 & 1 & 1 & 1 & 1000 & 1000 & 1 & 1 & 0.9 & 0.1 \\
2b & 1 & 1 & 1 & 1 & 1000 & 1000 & 1 & 1 & 0.9 & 0.1 \\
3 & 719.68 & 719.68 & $-1800$ & 1800 & 1320 & 1320 & $-1800$ & 1800 & 0.9 & 0.9 \\
4 & 138.94 & 138.94 & 3000 & $-3000$ & 1040.4 & 1040.4 & 3000 & $-3000$ & 0.9 & 0.9 \\
5 & 26.82 & 1 & 100 & 0 & 1004 & 1000 & 100 & 0 & 0.85 & 0.95 \\
6 & 1.01 & 1 & 25 & 25 & 1010 & 1000 & 25 & 25 & 0.8 & 0.9 \\
7 & 1.03546 & 1 & 0 & 0 & $1000 + 3.57 \cdot 10^{-8}$ & 1000 & 0 & 0 & 0.5 & 0.5 \\
\hline \\
\end{tabular}
\caption{Initial conditions for the Tests. All variables are reported in standard SI units: density $\rho$ in $\text{kg/m}^3$ and velocity $u$ in $\text{m/s}$, while the phase volume fraction $\alpha$ is dimensionless.}
\label{table:initial_condition}
\end{table}

\begin{table}[ht]
\centering
\footnotesize
\begin{tabular}{lccccccccc}
\hline
Test & $\rho_{01}$ & $\rho_{02}$ & $C_1$ & $C_2$ & $p_0$ & $\gamma_1$ & $\gamma_2$ & $T_{\text{end}}$ & L \\
\hline
1 & 1 & 1 & 1.18 & 1.41 & 0.995 & 1.4 & 2.0 & 0.125 & 1 \\
2a & 1 & $10^3$ & 1 & 5 & 0.714 & 1.4 & 2.8 & 0.05 & 1 \\
2b & 1 & $10^3$ & $10^2$ & $10^3$ & $7.1428 \cdot 10^3$ & 1.4 & 2.8 & 0.001 & 1 \\
3 & 1 & $10^3$ & $3.7416 \cdot 10^2$ & $1.5427 \cdot 10^3$ & $1 \cdot 10^5$ & 1.4 & 2.8 & $1.1 \cdot 10^{-5}$ & 0.1 \\
4 & 1 & $10^3$ & $3.7416 \cdot 10^2$ & $1.5427 \cdot 10^3$ & $1 \cdot 10^5$ & 1.4 & 2.8 & $2.2 \cdot 10^{-4}$ & 1 \\
5 & 1 & $10^3$ & $3.7416 \cdot 10^2$ & $1.534 \cdot 10^3$ & $1 \cdot 10^5$ & 1.4 & 2.8 & $2 \cdot 10^{-4}$ & 1 \\
6 & 1 & $10^3$ & 100 & $1 \cdot 10^3$ & $7.1428 \cdot 10^3$ & 1.4 & 2.8 & $1 \cdot 10^{-3}$ & 5 \\
7 & 1 & $10^3$ & 10 & $10^4$ & 71.43 & 1.4 & 2.8 & $2 \cdot 10^{-3}$ & 1 \\
\hline \\
\end{tabular}
\caption{Reference quantities and domain parameters. Units are specified in standard SI: density in $\text{kg/m}^3$, sound speed in $\text{m/s}$, pressure in $\text{Pa}$, final time $T_{\text{end}}$ in $\text{s}$, and the domain length $L$ in $\text{m}$.}
\label{table:ref_quantities}
\end{table}

We first consider a set of tests that serve to validate the schemes, in which the maximal wave speed of the full system, $\lambda_{\text{max}}^{\text{full}}$, and that of the split systems, $\lambda_{\text{max}}^{\text{split}}$, are of the same order of magnitude, so that the time steps given by the respective CFL conditions are very close. Then, we consider tests for which $\lambda_{\text{max}}^{\text{full}}$ is significantly larger than $\lambda_{\text{max}}^{\text{split}}$, so that the semi-implicit schemes are subject to a less restrictive CFL condition and can use larger time steps. For the Abgrall condition, both regimes are considered together in Test 2.

In all tests, we set $N_x = 500$ and $\nu_{\text{CFL}} = 0.6$, except for Test 6, which is run on a larger domain with $N_x = 2500$ cells, see Section~\ref{test_55}, and we compare the results of the two semi-implicit schemes with those of the explicit Rusanov scheme on the same grid and with the exact solution (Tests 2a, 2b, 3, 4 and 7) or with a numerical reference solution computed with the explicit Rusanov scheme on a fine grid of $100000$ cells (Tests 1, 5 and 6).
In all tests, homogeneous Neumann boundary conditions are imposed, both in the explicit and in the implicit stage. We do not consider relaxation and friction source terms, so that only the homogeneous part of system~\eqref{system:two_phase} is solved.

\subsection{Test 1: Generic Riemann Problem} \label{test_11}
\begin{figure}[htbp]
    \centering
    \includegraphics[width=\textwidth]{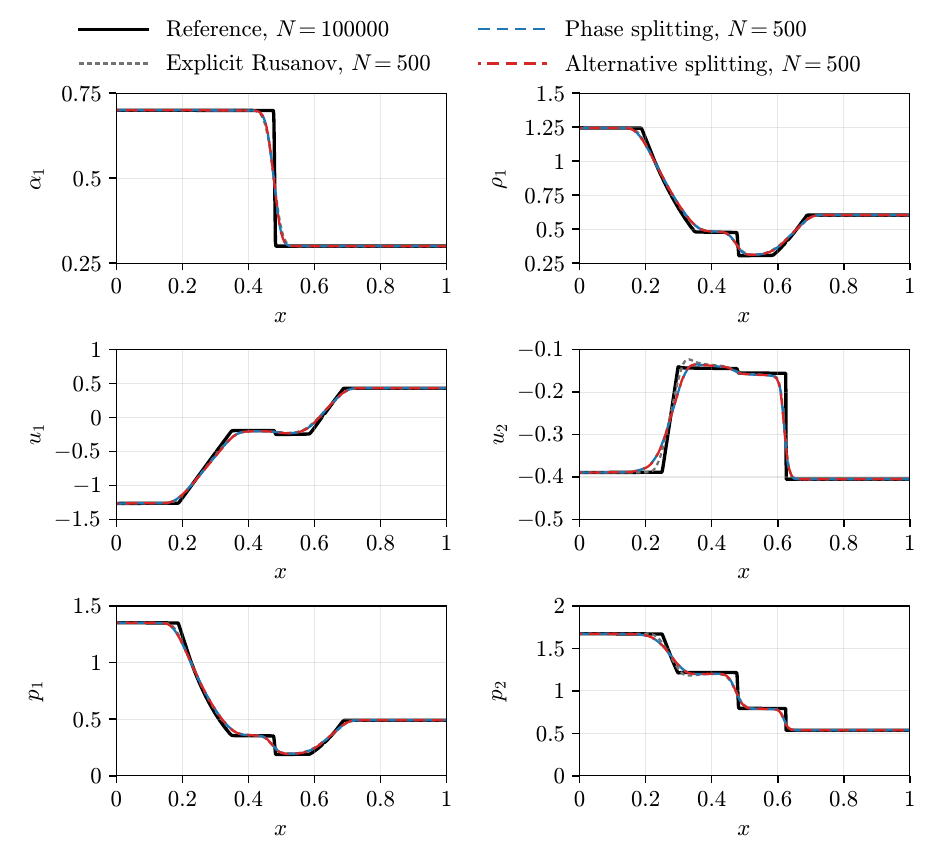}
    \caption{Numerical approximation of the two-phase flow model with the initial condition of Test 1 in Table \ref{table:initial_condition} at final time $T=0.125$.}
    \label{fig:test_11}
\end{figure}

This test consists of a generic Riemann problem with the initial discontinuity at the center of the domain ($x=0.5L$). It is a benchmark to assess the capability of the schemes to capture shocks, rarefactions and contact discontinuities in both phases in the absence of stiff terms. Since in this test the time step of the semi-implicit schemes coincides with the explicit one, see Table~\ref{table:time_step}, the three schemes yield practically the same results in Figure~\ref{fig:test_11}. They reproduce all waves of the reference solution at the correct positions, and the differences with respect to the reference solution are confined to the smearing of the discontinuities that is typical of first-order schemes on $500$ cells. 

\subsection{Test 2: Preservation of the Abgrall Condition} \label{test_1111}
Test 2 consists of two variants with identical initial data and different reference quantities, which assess the preservation of the pressure--velocity equilibrium in two acoustic regimes.

\subsubsection{Test 2a}
\begin{figure}[htbp]
    \centering
    \includegraphics[width=\textwidth]{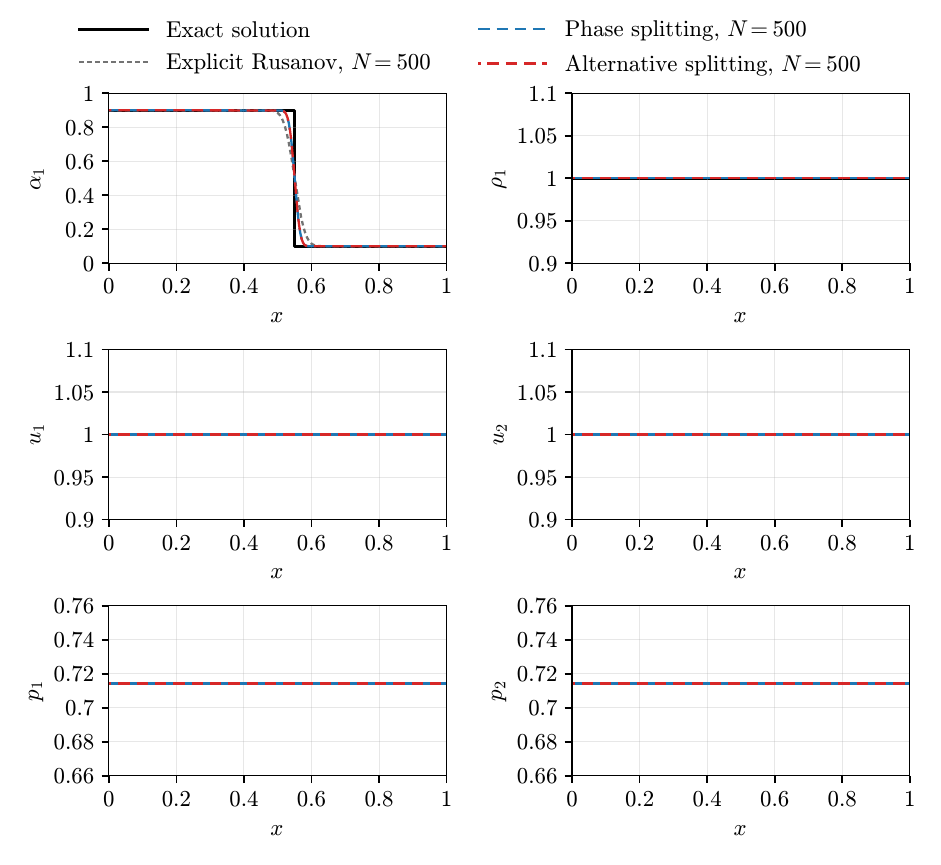}
    \caption{Numerical approximation of the two-phase flow model with the initial condition of Test 2a in Table \ref{table:initial_condition} at final time $T=0.05$.}
    \label{fig:test_1111}
\end{figure}

In this test, a jump in the volume fraction is advected with the uniform velocity $u_1=u_2=1$ and the uniform pressure $p_1=p_2=\rho_{01}C_1^2/\gamma_1\approx 0.714$ of both phases. This configuration is designed to assess whether the numerical schemes preserve the pressure--velocity equilibrium (\textit{Abgrall condition} \eqref{eq:Abgrall_condition}) also with the splitting approach. As shown in Figure \ref{fig:test_1111}, the splitting approach and the semi-implicit spatial discretization generate spurious oscillations neither in the pressure nor in the velocity of the two phases. The jump of the volume fraction is smeared over about $16$ cells by the semi-implicit schemes and over about $32$ cells by the explicit scheme. This is because the numerical viscosity of the explicit Rusanov scheme is computed from the largest wave speed of the full system, here the one of the liquid, whereas the explicit part of the splitting schemes does not involve the sound speed of the liquid.

\subsubsection{Test 2b}
\label{test_1104}
\begin{figure}[htbp]
    \centering
    \includegraphics[width=\textwidth]{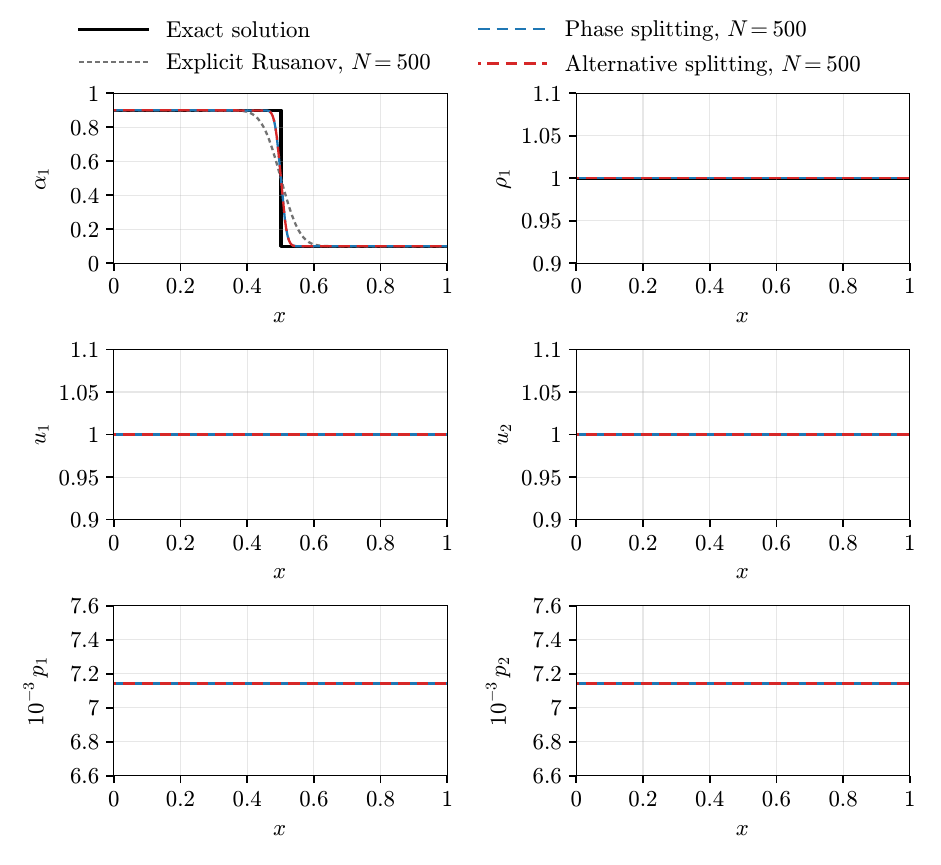}
    \caption{Numerical approximation of the two-phase flow model with the initial condition of Test 2b in Table \ref{table:initial_condition} at final time $T=0.001$.}
    \label{fig:test_1104}
\end{figure}

Test 2b uses the same initial data as Test 2a, with the density ratio $\rho_2/\rho_1 = 10^3$, but larger reference sound speeds and a larger reference pressure, see Table \ref{table:ref_quantities}. As shown in Figure \ref{fig:test_1104}, the splitting schemes preserve the uniform state up to round-off errors also in this case, \textit{i.e.} they remain ``well-balanced'' across different acoustic regimes. At the final time, the velocities of both phases computed with the two semi-implicit schemes differ from the uniform velocity by at most $10^{-14}$ in Test 2a and $3 \cdot 10^{-12}$ in Test 2b, and the pressures differ from the uniform pressure by at most $6 \cdot 10^{-11}$ and $2 \cdot 10^{-10}$ in relative terms, respectively. Here the jump of the volume fraction is smeared over about $19$ cells by the semi-implicit schemes and over about $57$ cells by the explicit scheme.

\subsection{Test 3: Double Rarefaction} \label{test_1}
\begin{figure}[htbp]
    \centering
    \includegraphics[width=\textwidth]{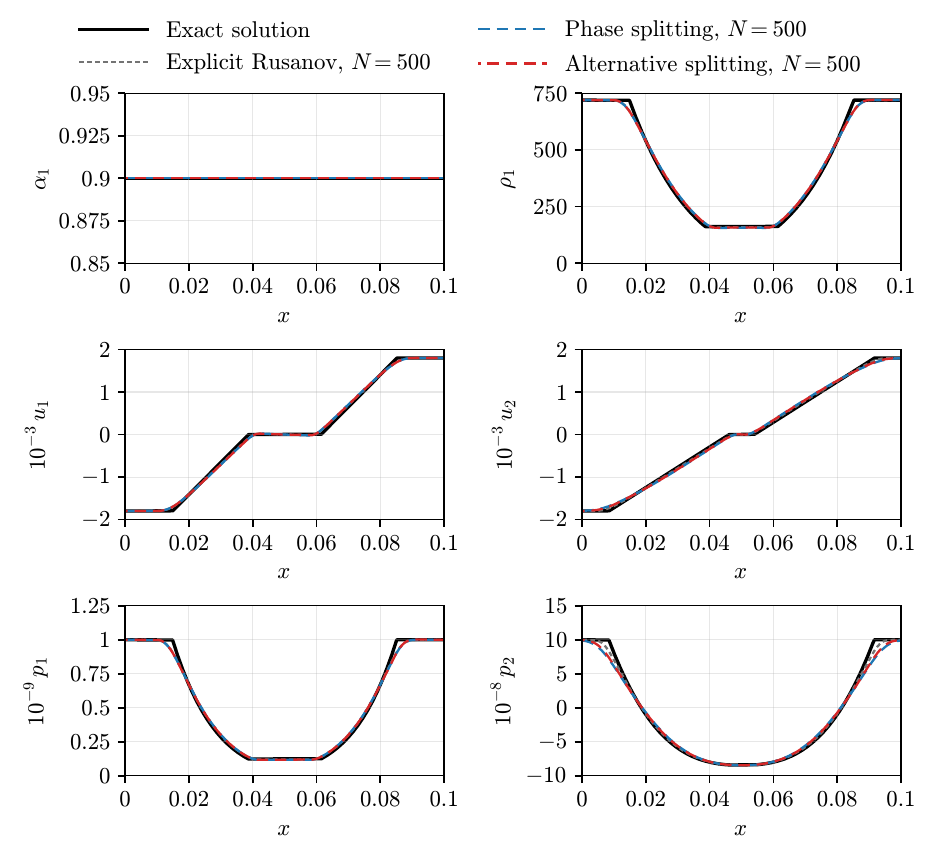}
    \caption{Numerical approximation of the two-phase flow model with the initial condition of Test 3 in Table \ref{table:initial_condition} at final time $T=1.1 \cdot 10^{-5}$.}
    \label{fig:test_1}
\end{figure}

This test is the four-rarefaction test proposed in \cite{romenski2004}, in which the reference densities differ by three orders of magnitude and the reference sound speeds by a factor of about four. Although the parameters are those of air and water, the initial states are extreme (gas at $\rho_1 \approx 720$ kg/m$^3$ and $p_1 = 10^9$ Pa, velocities $\pm 1800$ m/s), so that this is a severe robustness test. In the center, the liquid reaches $\rho_2 = 198.9$~kg/m$^3$ and $p_2 = -8.41\cdot 10^8$~Pa, \textit{i.e.} the vicinity of the lower limit $-p_\infty = -8.50\cdot 10^8$~Pa of the stiffened gas law \eqref{eq:eos_p2} for $\rho_2 \to 0$. We impose different initial velocity fields ($u_L <0 $ and $u_R > 0$) symmetrically in both phases. Since $\alpha_1$ is constant and only rarefaction waves occur, the model decouples into two barotropic Euler systems, one for each phase, and the exact solution is the superposition of their exact Riemann solutions \cite{theindumbser2022}. All three schemes agree well with the exact solution, see Figure~\ref{fig:test_1}. The semi-implicit schemes are only slightly more diffusive at the heads of the liquid rarefactions, since the liquid acoustic waves are integrated implicitly.

\subsection{Test 4: Double Shock} 
\label{test_2}

\begin{figure}[htbp]
    \centering
    \includegraphics[width=\textwidth]{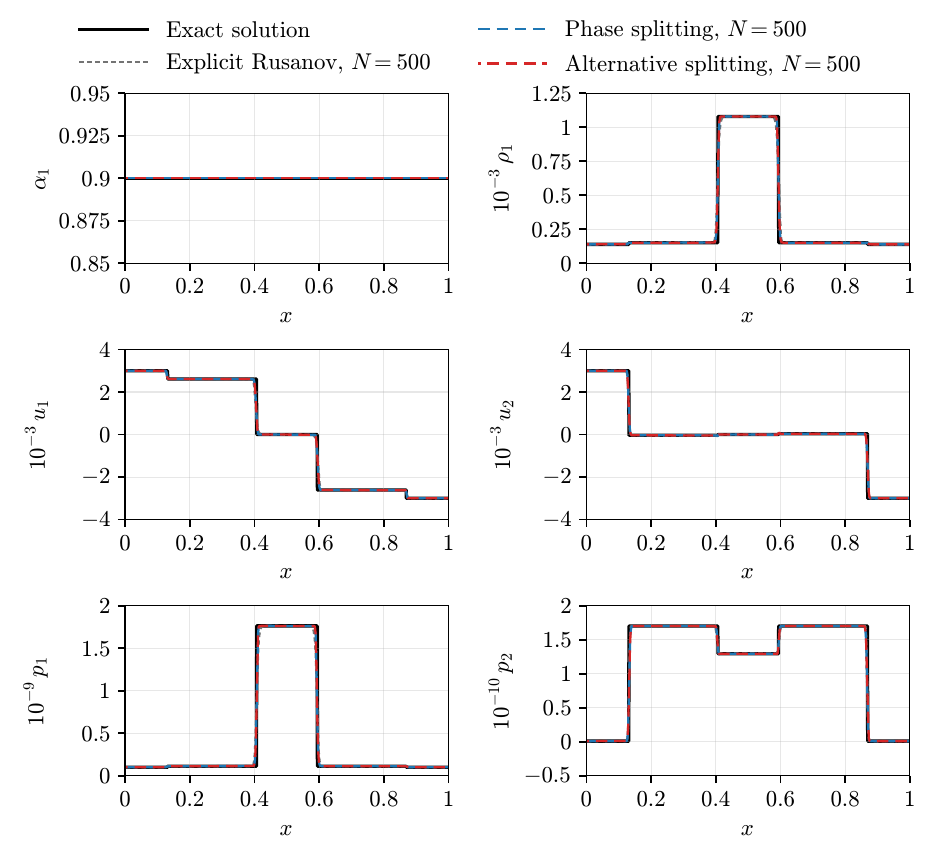}
    \caption{Numerical approximation of the two-phase flow model with the initial condition of Test 4 in Table \ref{table:initial_condition} at final time $T=2.2 \cdot 10^{-4}$.}
    \label{fig:test_2}
\end{figure}

In contrast to the previous configuration in Subsection \ref{test_1}, we consider here high-speed converging flows ($u_L > 0 $ and $u_R < 0$), as proposed in \cite{romenski2004}. In Figure \ref{fig:test_2} we observe that two strong propagating shock waves are produced and the density moves towards the domain center. The volume fraction remains constant, $\alpha_1 = 0.9$, since it is only transported. Each phase produces two shocks: the fast liquid shocks bring the liquid almost to rest, while the slow gas shocks stop and strongly compress the gas in the center of the domain. Across the gas shocks, the Rankine--Hugoniot condition of the relative velocity equation \eqref{eq:rel_vel}, see \cite{theindumbser2022}, couples the deceleration and compression of the gas to the enthalpy of the liquid, so that the liquid density and pressure decrease in the center, although the mixture as a whole is compressed. We compute the exact solution with the exact Riemann solver of \cite{theindumbser2022}. All three schemes capture the four shocks at the positions of the exact solution. The gas shocks are resolved more sharply by the semi-implicit schemes, since the numerical viscosity of their explicit part does not depend on the sound speed of the liquid, whereas the liquid shocks are captured by all three schemes with practically the same accuracy.

%\subsection{Group 2}

\subsection{Test 5: Strong Shock Tube} \label{test_33}

\begin{figure}[htbp]
    \centering
    \includegraphics[width=\textwidth]{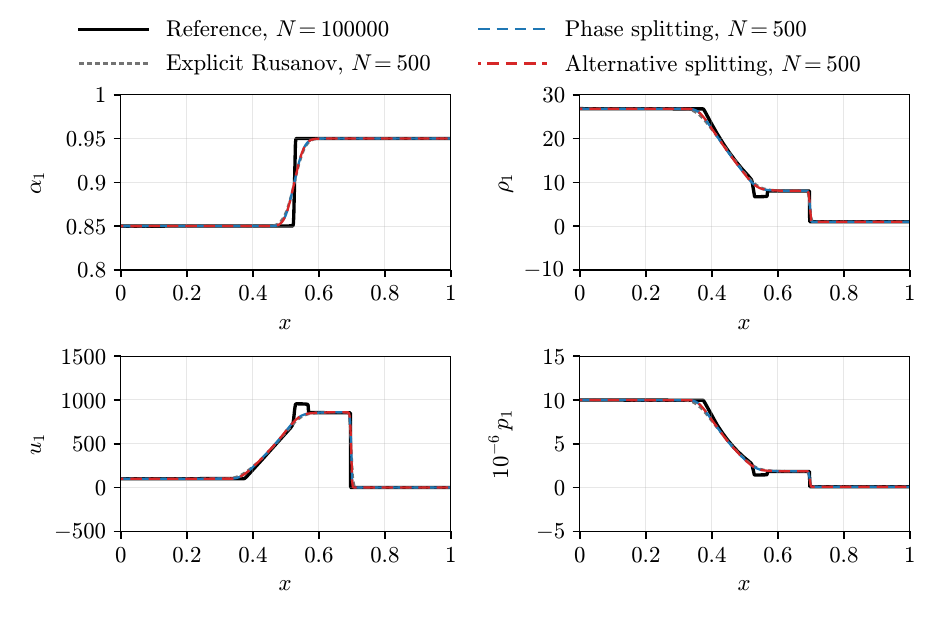}
    \caption{Numerical approximation of the two-phase flow model with the initial condition of Test 5 in Table \ref{table:initial_condition} at final time $T=2 \cdot 10^{-4}$.}
    \label{fig:test_33}
\end{figure}

This test considers a two-phase shock tube with an expanding left state characterized by high pressure and high density. The right state is at rest. Figure \ref{fig:test_33} shows how the schemes handle the transient and capture the propagation speed of the resulting shock in the gas--liquid mixture. All three schemes capture the positions of the waves of the reference solution. The gas waves are resolved somewhat more sharply by the semi-implicit schemes. None of the schemes resolves on $500$ cells the narrow region of width about $0.04$ to the right of the material interface at $x \approx 0.53$, in which the gas reaches its maximum velocity $u_1 \approx 950$ in the reference solution.

\subsection{Test 6: Moving Acoustic Wave} \label{test_55}

\begin{figure}[htbp]
    \centering
    \includegraphics[width=0.8\textwidth]{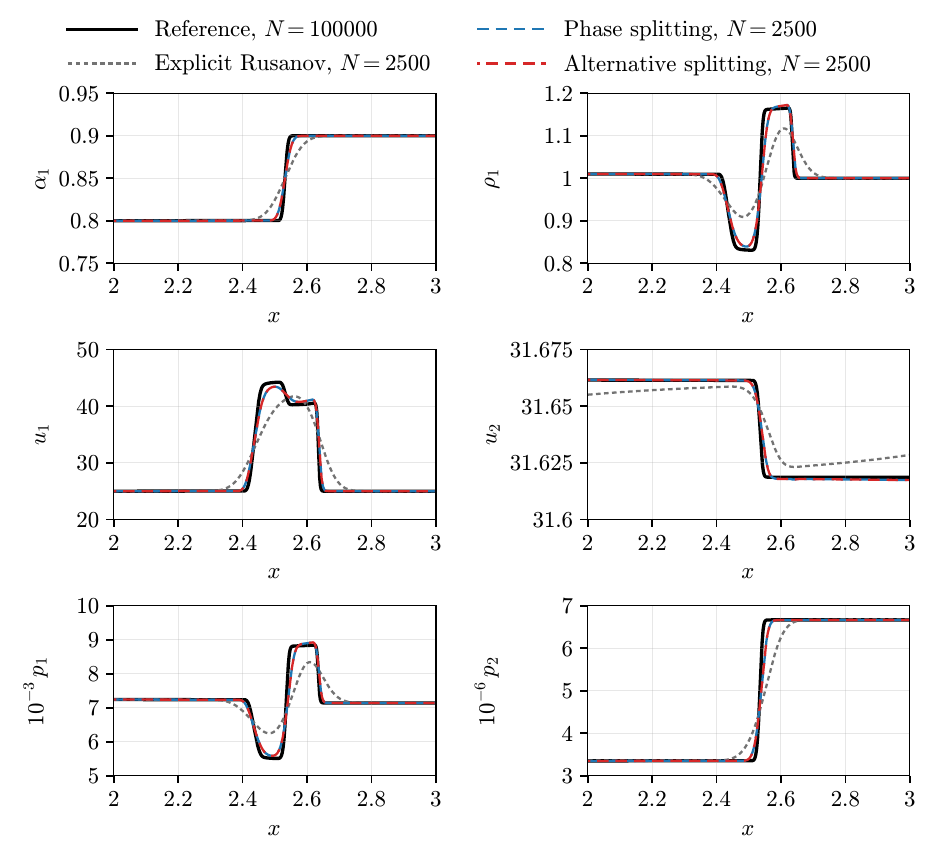}
    \caption{Numerical approximation of the two-phase flow model with the initial condition of Test 6 in Table \ref{table:initial_condition} at final time $T=1 \cdot 10^{-3}$, shown on the subdomain $[2,3]$ of the computational domain $[0,5]$.}
    \label{fig:test_55}
\end{figure}

This test is designed to assess the capability of the schemes to capture the acoustic waves of the gas phase. Here, the time step $\Delta t$ is dictated only by the explicit CFL condition based on the gas and material velocities ($|u| + c_1$). Consequently, the slow material waves and the gas acoustic waves are well resolved and accurately preserved. In contrast, the fast acoustic waves of the liquid phase travel at a much higher speed ($|u| + c_2$) and correspond to a local Courant number well above unity. Consequently, the under-resolved fast acoustic modes are strongly damped by the implicit time integration with large time steps. In order to prevent these waves from reaching the boundaries before the final time, this test is run on the larger domain $[0,5]$ with $N_x = 2500$ cells, \textit{i.e.} with the same mesh spacing as the other tests, and Figure~\ref{fig:test_55} shows only the subdomain $[2,3]$. The volume fraction and the gas quantities computed with the two semi-implicit schemes are much closer to the reference solution than those of the explicit scheme, which smears the material interface and the gas waves, since its numerical viscosity is computed from the sound speed of the liquid. The liquid states of both semi-implicit schemes agree with the reference solution.

\subsection{Test 7: Large Ratio of the Sound Speeds}
This test illustrates that the semi-implicit splitting schemes are not limited by a CFL condition based on the sound speed of the phase that is treated implicitly, here the liquid, and are therefore suitable for low Mach numbers in that phase. We consider a mixture with equal volume fractions, $\alpha_1=\alpha_2=0.5$, and a low-Mach acoustic pressure perturbation in a fluid at rest. To maximize the disparity of the scales, we couple a very stiff liquid phase ($c_2=10^4$) with a slow gas phase ($c_1=10$). The initial data are in pressure equilibrium, with $p_1 = p_2 = 1.05\,p_0$ for $x < L/2$ and $p_1 = p_2 = p_0$ for $x \geq L/2$.

In a standard fully explicit framework, the global time step that ensures stability is governed by the fastest wave in the domain, \textit{i.e.} by the very large sound speed of the liquid. This imposes a prohibitive restriction, $\Delta t \approx 10^{-7}$, whereas in the semi-implicit approach the stable time step is dictated by the gas phase, $\Delta t \approx 10^{-4}$. As a result, the semi-implicit schemes capture the long-term propagation of the gas acoustic waves with time steps much larger than the explicit limit.

\begin{figure}[htbp]
    \centering
    \includegraphics[width=\textwidth]{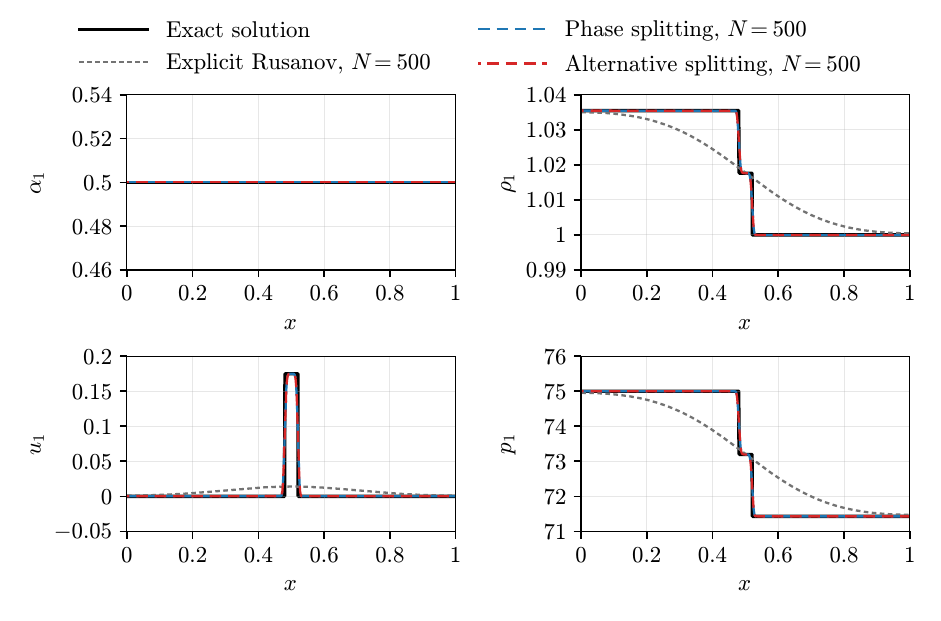}
    \caption{Numerical approximation of the two-phase flow model with the initial condition of Test 7 in Table \ref{table:initial_condition} at final time $T=2 \cdot 10^{-3}$.}
    \label{fig:test_62}
\end{figure}

Figure \ref{fig:test_62} shows that the semi-implicit schemes resolve the gasdynamics much better than the explicit Rusanov scheme. Indeed, the Rusanov scheme computes the numerical viscosity from the liquid sound speed, which introduces excessive dissipation and smears the gas waves. We compute the exact solution with the exact Riemann solver of \cite{theindumbser2022}. It consists of a rarefaction and a shock in each phase, and at the final time the two liquid waves have left the domain. Between the two gas waves the exact solution has the gas velocity $u_1 = 0.175$, which the semi-implicit schemes reproduce on $500$ cells, whereas the explicit Rusanov scheme reaches only $u_1 \approx 0.014$ on $500$ cells and $u_1 \approx 0.08$ even on $20000$ cells. On the other hand, as expected, the semi-implicit schemes do not capture the fast liquid waves, which are damped by the large time step dictated by the gas phase.

\begin{table}[ht]
\centering
\footnotesize
\begin{tabular}{lcccc}
\hline
Test & $\Delta t_{\text{expl}}$ & $\Delta t_{\text{phase}}$ & $\Delta t_{\text{altern}}$ & $\Delta t_{\text{phase}} / \Delta t_{\text{expl}}$ \\
\hline
1 & $4.81 \cdot 10^{-4}$ & $4.81 \cdot 10^{-4}$ & $4.81 \cdot 10^{-4}$ & $1.00$\\
2a & $2.00 \cdot 10^{-4}$ & $6.00 \cdot 10^{-4}$ & $6.00 \cdot 10^{-4}$ & $3.00$\\
2b & $1.20 \cdot 10^{-6}$ & $1.19 \cdot 10^{-5}$ & $1.19 \cdot 10^{-5}$ & $9.92$\\
3 & $3.17 \cdot 10^{-8}$ & $3.76 \cdot 10^{-8}$ & $3.76 \cdot 10^{-8}$ & $1.19$\\
4 & $2.61 \cdot 10^{-7}$ & $3.00 \cdot 10^{-7}$ & $2.75 \cdot 10^{-7}$ & $1.15$\\
5 & $7.32 \cdot 10^{-7}$ & $1.46 \cdot 10^{-6}$ & $1.46 \cdot 10^{-6}$ & $1.99$\\
6 & $1.16 \cdot 10^{-6}$ & $9.58 \cdot 10^{-6}$ & $9.58 \cdot 10^{-6}$ & $8.26$\\
7 & $1.20 \cdot 10^{-7}$ & $1.19 \cdot 10^{-4}$ & $1.19 \cdot 10^{-4}$ & $991.67$\\
\hline \\
\end{tabular}
\caption{Maximum time step over the whole simulation for the explicit scheme and for the two semi-implicit schemes, and ratio of the semi-implicit to the explicit time step. In all tests except Test 4, where $|\lambda_{1,4}|$ exceeds $|u_1|+c_1$, $\Delta t_{\text{altern}} = \Delta t_{\text{phase}}$. In all simulations, $\nu_{\text{CFL}} = 0.6$ and $N_x = 500$, except for Test 6 with $N_x = 2500$ on the larger domain.}
\label{table:time_step}
\end{table}

\clearpage
\section{Conclusions}
\label{Section:conclusions}
In this paper, we have proposed two different flux vector splittings and two semi-implicit finite volume schemes for the one-dimensional hyperbolic and thermodynamically compatible (HTC) isentropic two-\hspace{0pt}phase flow model of Romenski \textit{et al.} \cite{romenski2004, romenski2007, toro2010}, for gas--liquid mixtures in which the liquid is weakly compressible. Both schemes are based on a splitting of the flux at the continuous level. In the phase splitting, the gas-system is discretized explicitly and the entire liquid-system implicitly, which leads to a nonlinear equation system for the liquid pressure that also contains the convective terms of the liquid. In the alternative splitting, the convective terms of the liquid are moved into the explicit stage, so that only the acoustic waves of the liquid are treated implicitly, and the corresponding equation system for the pressure can be solved without nested Picard iterations.

We have analysed the eigenstructure of the resulting subsystems. The explicit subsystems are strongly hyperbolic except for a few degenerate states, while the implicit subsystems are only weakly hyperbolic for a nonzero liquid pressure, a property that they share with the explicit convective subsystem of the flux vector splitting of Toro and V{\'a}zquez-Cend{\'o}n \cite{torovazquez2012}. For both schemes, we have proven that the pressure--velocity equilibrium across material interfaces (\textit{Abgrall condition}) is preserved exactly, which is also confirmed by the numerical tests. Moreover, for well-prepared initial data, both schemes are consistent with the volume balance of the low Mach number limit of the liquid phase.

The numerical tests show that the semi-implicit part of the schemes is not limited by a CFL condition based on the sound speed of the liquid, so that the time step is dictated only by the gas waves, by the material waves and, in the alternative splitting, by the convective velocity of the liquid bulk, and the schemes are therefore suitable for low Mach numbers in the liquid phase. The gain in the time step grows as the Mach number of the liquid phase decreases, and the semi-implicit time step can be up to three to four orders of magnitude larger than the prohibitively small explicit one. At the same time, the gas waves are resolved much better than with the explicit Rusanov scheme, whose numerical viscosity is dictated by the sound speed of the liquid and leads to an excessive dissipation of the slow waves. The price to pay is that the fast acoustic waves of the liquid are damped by the large time steps, which is acceptable whenever the main interest lies in the gasdynamics and in the material waves.

The schemes presented in this paper are first-order accurate and one-dimensional. Future work will concern their extension to higher order of accuracy in space and time and to multiple space dimensions.

\section*{Acknowledgement}
\noindent TT would like to thank Sara Rinaldi for fruitful discussions on this topic.\\ GP and TT received funding from the European Union's Horizon Europe research and innovation program under the Marie Sk{\l}odowska-Curie Doctoral Network DataHyking (Grant No. 101072546). MD, GP and AS are supported by the Ministry of University
and Research under the PRIN-2022 project on ``High Order
Structure-preserving Semi-implicit Schemes for Hyperbolic
Equations'' (Project No. 2022JH87B4). MD was funded by the Fondazione Caritro via the project SOPHOS and by the European Research Council (ERC) under the European Union's Horizon Europe research and innovation programme via the project SOPHOS, grant agreement no. ERC-ADG-2025-101265878-SOPHOS. Views and opinions expressed are however those of the authors only and do not necessarily reflect those of the European Union or the European Research Council Executive Agency. Neither the European Union nor the granting authority can be held responsible for them. MD, GP and TT are members of the INdAM Research National Group of Scientific Computing (INdAM-GNCS).

\subsection*{Data Availability} The numerical data supporting the results of this study are available from the corresponding author upon reasonable request.

\subsection*{Conflict of interest} The authors declare no conflict of interest.

\section*{Dedication}
\textbf{This paper is dedicated to Prof. Dr. Christian Klingenberg in honor of his 70\textsuperscript{th} birthday.}

\bibliography{References}

@article{puppo2023,
  author  = {M. Luk{\'a}{\v{c}}ov{\'a}-Medvi{\v{d}}ov{\'a} and G. Puppo and A. Thomann},
  title   = {An all {M}ach number finite volume method for isentropic two-phase flow},
  journal = {Journal of Numerical Mathematics},
  volume  = {31},
  pages   = {175--204},
  year    = {2023}
}

@article{thomann2023,
  author  = {A. Thomann and M. Dumbser},
  title   = {Thermodynamically compatible discretization of a compressible two-fluid model with two entropy inequalities},
  journal = {Journal of Scientific Computing},
  volume  = {97},
  pages   = {9},
  year    = {2023}
}

@article{theindumbser2022,
  author  = {F. Thein and E. Romenski and M. Dumbser},
  title   = {Exact and numerical solutions of the {R}iemann problem for a conservative model of compressible two-phase flows},
  journal = {Journal of Scientific Computing},
  volume  = {93},
  pages   = {83},
  year    = {2022}
}

@article{dumbser2011,
  author  = {M. Dumbser},
  title   = {A simple two-phase method for the simulation of complex free surface flows},
  journal = {Computer Methods in Applied Mechanics and Engineering},
  volume  = {200},
  pages   = {1204--1219},
  year    = {2011}
}

@article{toro2010,
  author  = {E. Romenski and D. Drikakis and E. F. Toro},
  title   = {Conservative models and numerical methods for compressible two-phase flow},
  journal = {Journal of Scientific Computing},
  volume  = {42},
  pages   = {68--95},
  year    = {2010}
}

@article{saurel1999,
  author  = {R. Saurel and R. Abgrall},
  title   = {A multiphase {G}odunov method for compressible multifluid and multiphase flows},
  journal = {Journal of Computational Physics},
  volume  = {150},
  pages   = {425--467},
  year    = {1999}
}

@article{abgrall1996,
  author  = {R. Abgrall},
  title   = {How to prevent pressure oscillations in multicomponent flow calculations: a quasi conservative approach},
  journal = {Journal of Computational Physics},
  volume  = {125},
  pages   = {150--160},
  year    = {1996}
}

@article{dumbser2010,
  author  = {M. Dumbser and A. Hidalgo and M. Castro and C. Par{\'e}s and E. F. Toro},
  title   = {{FORCE} schemes on unstructured meshes {II}: non-conservative hyperbolic systems},
  journal = {Computer Methods in Applied Mechanics and Engineering},
  volume  = {199},
  pages   = {625--647},
  year    = {2010}
}

@article{romenski2004,
  author  = {E. Romenski and E. F. Toro},
  title   = {Compressible two-phase flow models: two-pressure models and numerical methods},
  journal = {Computational Fluid Dynamics Journal},
  volume  = {13},
  pages   = {403--416},
  year    = {2004}
}

@article{romenski2007,
  author  = {E. Romenski and A. D. Resnyansky and E. F. Toro},
  title   = {Conservative hyperbolic formulation for compressible two-phase flow with different phase pressures and temperatures},
  journal = {Quarterly of Applied Mathematics},
  volume  = {65},
  pages   = {259--279},
  year    = {2007}
}

@article{puppo2026,
  author  = {G. Puppo and T. Rey and T. Tenna},
  title   = {Formal derivation of an isentropic two-phase flow model from the multispecies {B}oltzmann equation},
  journal = {Physical Review E},
  volume  = {113},
  pages   = {034119},
  year    = {2026}
}

@article{flatten2011,
  author  = {T. Fl{\aa}tten and H. Lund},
  title   = {Relaxation two-phase flow models and the subcharacteristic condition},
  journal = {Mathematical Models and Methods in Applied Sciences},
  volume  = {21},
  pages   = {2379--2407},
  year    = {2011}
}

@article{baer1986,
  author  = {M. R. Baer and J. W. Nunziato},
  title   = {A two-phase mixture theory for the deflagration-to-detonation transition ({DDT}) in reactive granular materials},
  journal = {International Journal of Multiphase Flow},
  volume  = {12},
  pages   = {861--889},
  year    = {1986}
}

@article{lund2012,
  author  = {H. Lund},
  title   = {A hierarchy of relaxation models for two-phase flow},
  journal = {SIAM Journal on Applied Mathematics},
  volume  = {72},
  pages   = {1713--1741},
  year    = {2012}
}

@incollection{godunov1995,
  author    = {S. K. Godunov and E. I. Romenski},
  title     = {Thermodynamics, conservation laws, and symmetric forms of differential equations in mechanics of continuous media},
  booktitle = {Computational Fluid Dynamics Review 95},
  publisher = {John Wiley},
  address   = {New York},
  pages     = {19--31},
  year      = {1995}
}

@article{saurel2001,
  author  = {R. Saurel and O. Lemetayer},
  title   = {A multiphase model for compressible flows with interfaces, shocks, detonation waves and cavitation},
  journal = {Journal of Fluid Mechanics},
  volume  = {431},
  pages   = {239--271},
  year    = {2001}
}

@techreport{harlow1970,
  author      = {F. H. Harlow and A. A. Amsden},
  title       = {Fluid dynamics. {A} {LASL} monograph},
  institution = {Los Alamos Scientific Laboratory},
  number      = {LA-4700},
  year        = {1971}
}

@article{delrio2025,
  author  = {L. R{\'{\i}}o-Mart{\'{\i}}n and F. Dhaouadi and M. Dumbser},
  title   = {An exactly curl-free finite-volume/finite-difference scheme for a hyperbolic compressible isentropic two-phase model},
  journal = {Journal of Scientific Computing},
  volume  = {102},
  pages   = {13},
  year    = {2025}
}

@article{ferrari2024,
  author  = {D. Ferrari and M. Dumbser},
  title   = {A semi-implicit finite volume scheme for incompressible two-phase flows},
  journal = {Communications on Applied Mathematics and Computation},
  volume  = {6},
  pages   = {2295--2330},
  year    = {2024}
}

@article{herard2019,
  author  = {J.-M. H{\'e}rard and H. Mathis},
  title   = {A three-phase flow model with two miscible phases},
  journal = {ESAIM: Mathematical Modelling and Numerical Analysis},
  volume  = {53},
  pages   = {1373--1389},
  year    = {2019}
}

@article{gallouet2004,
  author  = {T. Gallou{\"e}t and J.-M. H{\'e}rard and N. Seguin},
  title   = {Numerical modeling of two-phase flows using the two-fluid two-pressure approach},
  journal = {Mathematical Models and Methods in Applied Sciences},
  volume  = {14},
  pages   = {663--700},
  year    = {2004}
}

@article{coquel2002,
  author  = {F. Coquel and T. Gallou{\"e}t and J.-M. H{\'e}rard and N. Seguin},
  title   = {Closure laws for a two-fluid two-pressure model},
  journal = {Comptes Rendus Math{\'e}matique},
  volume  = {334},
  pages   = {927--932},
  year    = {2002}
}

@article{ferrari2025,
  author  = {D. Ferrari and I. Peshkov and E. Romenski and M. Dumbser},
  title   = {A unified {HTC} multiphase model of continuum mechanics},
  journal = {Journal of Computational Physics},
  volume  = {521},
  pages   = {113553},
  year    = {2025}
}

@article{ferrari2026,
  author  = {D. Ferrari and I. Peshkov and M. Dumbser},
  title   = {A configurational numerical viscosity for thermodynamically compatible finite volume schemes of compressible multiphase flows},
  journal = {arXiv preprint arXiv:2610.23247},
  year    = {2026}
}

@article{saurel2018,
  author  = {R. Saurel and C. Pantano},
  title   = {Diffuse-interface capturing methods for compressible two-phase flows},
  journal = {Annual Review of Fluid Mechanics},
  volume  = {50},
  pages   = {105--130},
  year    = {2018}
}

@book{ishii1975,
  author    = {M. Ishii},
  title     = {Thermo-Fluid Dynamic Theory of Two-Phase Flow},
  series    = {Collection de la Direction des {\'E}tudes et Recherches d'{\'E}lectricit{\'e} de France},
  volume    = {22},
  publisher = {Eyrolles},
  address   = {Paris},
  year      = {1975}
}

@book{prosperetti2009,
  editor    = {A. Prosperetti and G. Tryggvason},
  title     = {Computational Methods for Multiphase Flow},
  publisher = {Cambridge University Press},
  address   = {Cambridge},
  year      = {2007}
}

@article{casulli1984,
  author  = {V. Casulli and D. Greenspan},
  title   = {Pressure method for the numerical solution of transient, compressible fluid flows},
  journal = {International Journal for Numerical Methods in Fluids},
  volume  = {4},
  pages   = {1001--1012},
  year    = {1984}
}

@article{casulli1990,
  author  = {V. Casulli},
  title   = {Semi-implicit finite difference methods for the two-dimensional shallow water equations},
  journal = {Journal of Computational Physics},
  volume  = {86},
  pages   = {56--74},
  year    = {1990}
}

@article{klein1995,
  author  = {R. Klein},
  title   = {Semi-implicit extension of a {G}odunov-type scheme based on low {M}ach number asymptotics {I}: one-dimensional flow},
  journal = {Journal of Computational Physics},
  volume  = {121},
  pages   = {213--237},
  year    = {1995}
}

@article{munz2003,
  author  = {C.-D. Munz and S. Roller and R. Klein and K. J. Geratz},
  title   = {The extension of incompressible flow solvers to the weakly compressible regime},
  journal = {Computers \& Fluids},
  volume  = {32},
  pages   = {173--196},
  year    = {2003}
}

@article{parkmunz2005,
  author  = {J. H. Park and C.-D. Munz},
  title   = {Multiple pressure variables methods for fluid flow at all {M}ach numbers},
  journal = {International Journal for Numerical Methods in Fluids},
  volume  = {49},
  pages   = {905--931},
  year    = {2005}
}

@article{brugnanocasulli2008,
  author  = {L. Brugnano and V. Casulli},
  title   = {Iterative solution of piecewise linear systems},
  journal = {SIAM Journal on Scientific Computing},
  volume  = {30},
  pages   = {463--472},
  year    = {2008}
}

@article{casullizanolli2012,
  author  = {V. Casulli and P. Zanolli},
  title   = {Iterative solutions of mildly nonlinear systems},
  journal = {Journal of Computational and Applied Mathematics},
  volume  = {236},
  pages   = {3937--3947},
  year    = {2012}
}

@article{dumbsercasulli2016,
  author  = {M. Dumbser and V. Casulli},
  title   = {A conservative, weakly nonlinear semi-implicit finite volume scheme for the compressible {N}avier--{S}tokes equations with general equation of state},
  journal = {Applied Mathematics and Computation},
  volume  = {272},
  pages   = {479--497},
  year    = {2016}
}

@article{torovazquez2012,
  author  = {E. F. Toro and M. E. V{\'a}zquez-Cend{\'o}n},
  title   = {Flux splitting schemes for the {E}uler equations},
  journal = {Computers \& Fluids},
  volume  = {70},
  pages   = {1--12},
  year    = {2012}
}

@article{pareschirusso2005,
  author  = {L. Pareschi and G. Russo},
  title   = {Implicit-explicit {R}unge--{K}utta schemes and applications to hyperbolic systems with relaxation},
  journal = {Journal of Scientific Computing},
  volume  = {25},
  pages   = {129--155},
  year    = {2005}
}

@article{thomann2020,
  author  = {A. Thomann and M. Zenk and G. Puppo and C. Klingenberg},
  title   = {An all speed second order {IMEX} relaxation scheme for the {E}uler equations},
  journal = {Communications in Computational Physics},
  volume  = {28},
  pages   = {591--620},
  year    = {2020}
}

@article{castro2006,
  author  = {M. J. Castro and J. M. Gallardo and C. Par{\'e}s},
  title   = {High-order finite volume schemes based on reconstruction of states for solving hyperbolic systems with nonconservative products. {A}pplications to shallow-water systems},
  journal = {Mathematics of Computation},
  volume  = {75},
  pages   = {1103--1134},
  year    = {2006}
}

@article{pares2006,
  author  = {C. Par{\'e}s},
  title   = {Numerical methods for nonconservative hyperbolic systems: a theoretical framework},
  journal = {SIAM Journal on Numerical Analysis},
  volume  = {44},
  pages   = {300--321},
  year    = {2006}
}
\bibliographystyle{acm}

\end{document}